\documentclass[DIV=13]{scrartcl}
\usepackage[sb]{libertine} %
\usepackage[T1]{fontenc}
\usepackage{textcomp}
\usepackage[varqu,varl]{zi4}%
\usepackage[amsthm,upint]{libertinust1math} %
\usepackage[scr=boondoxo,bb=boondox]{mathalpha} %
\usepackage{bm}
\usepackage{mathtools}
\usepackage{todonotes}
\usepackage{booktabs}
\usepackage{enumitem}
\usepackage{placeins}
\usepackage{colortbl}
\usepackage{bm}
\usepackage[style=numeric, giveninits]{biblatex}
\AtEveryBibitem{%
  \clearfield{urlyear}%
  \clearfield{urldate}%
}%
\AtEveryCitekey{%
  \clearfield{urlyear}%
  \clearfield{urldate}%
}%
\usepackage{xcolor}
\usepackage{float}
\RequirePackage{siunitx}
\RequirePackage[algoruled,norelsize, linesnumbered, lined, commentsnumbered]{algorithm2e}
\usepackage[format=plain,labelfont=bf]{caption}
\usepackage{subcaption}
\usepackage[pdftex,colorlinks=true,linkcolor={blue!50!black},citecolor={red!50!black},urlcolor=blue]{hyperref}

\newcommand{\symbf}[1]{\bm{#1}}
\newcommand{\symcal}[1]{\mathcal{#1}}
\newcommand{\symbfit}[1]{\boldsymbol{#1}}
\newcommand{\symbfcal}[1]{\boldsymbol{\mathcal #1}}
\newcommand\restrict[1]{\raisebox{-.5ex}{$|$}_{#1}}

\definecolor{myblue}{RGB}{0 83 139}
\definecolor{myred}{RGB}{114 16 69}
\definecolor{mygreen}{RGB}{0 94 0}
\usepackage{mycommands}

\usepackage{booktabs}
\author{Nils Margenberg
  \thanks{Helmut Schmidt University,
    Faculty of Mechanical and Civil Engineering,
    Holstenhofweg 85,
    22043 Hamburg,
    Germany,
    \href{mailto:margenbn@hsu-hh.de}{\texttt{margenbn@hsu-hh.de}} (Corresponding Author)}
  \and Peter Munch\thanks{Uppsala University,
    Department of Information Technology,
    Box 337, Uppsala, 75105, Sweden,
    \href{mailto:peter.munch@it.uu.se}{\texttt{peter.munch@it.uu.se}}}}
\date{}
\title{A Space-Time Multigrid Method for Space-Time Finite Element Discretizations of Parabolic and Hyperbolic PDEs}
\usepackage{biblatex}
\begin{document}
\maketitle
\vspace{-8ex}%
\begin{abstract}
We present a space-time multigrid method based on tensor-product space-time finite element discretizations. The method is facilitated by the matrix-free capabilities of the {\ttfamily deal.II} library. It addresses both high-order continuous and discontinuous variational time discretizations with spatial finite element discretizations. The effectiveness of multigrid methods in large-scale stationary problems is well established. However, their application in the space-time context poses significant challenges, mainly due to the construction of suitable smoothers. To address these challenges, we develop a space-time cell-wise additive Schwarz smoother and demonstrate its effectiveness on the heat and acoustic wave equations. 

The matrix-free framework of the {\ttfamily deal.II} library supports various multigrid strategies, including $h$-, $p$-, and $hp$-refinement across spatial and temporal dimensions. Extensive empirical evidence, provided through scaling and convergence tests on high-performance computing platforms, demonstrate high performance on perturbed meshes and problems with heterogeneous and discontinuous coefficients. Throughputs of over a billion degrees of freedom per second are achieved on problems with more than a trillion global degrees of freedom. The results prove that the space-time multigrid method can effectively solve complex problems in high-fidelity simulations and show great potential for use in coupled problems.

  \noindent\emph{MSC2020: 65M60, 65M55, 65F10, 65Y05}\\
  \emph{Keywords: space-time finite elements, space-time multigrid, matrix-free,
    high-order,\\
    high-performance computing, tensor-product}
\end{abstract}
\section{Introduction}
\label{sec:orgb89780c}
In the recent years, parallel computing has become increasingly important for
improving the efficiency and reducing computation times of numerical solutions
for partial differential equations (PDEs). One significant area of research has
been the development of space-time finite element methods and parallel time
integration techniques, which aim to take advantage of spatial and temporal
parallelism.

We present a space-time multigrid (STMG) method utilizing tensor-product
space-time finite elements within the {\ttfamily deal.II} library
framework~\cite{arndtDealIILibrary2023,kronbichlerGenericInterfaceParallel2012}.
This approach combines high-order variational time discretizations with spatial
finite elements to extend efficient multigrid methods to space-time problems.
Employing a space-time cell-wise additive Schwarz method (ASM) as the smoother,
we demonstrate high performance for heat and acoustic wave equations on
perturbed meshes and problems with heterogeneous coefficients. Our scalable
implementation is validated through tests on high-performance computing
platforms and available on GitHub~\url{https://github.com/nlsmrg/dealii-stfem}.

Further investigations on space-time finite element methods have been done
in~\cite{drflerSpaceTimeDiscontinuousGalerkin2016,drflerParallelAdaptiveDiscontinuous2019,langerAdaptiveSpaceTime2022},
where numerical results with an adaptive algorithm are presented. The advantages
of the variational time discretization are the natural integration with the
variational space discretization and the natural capture of coupled problems and
nonlinearities. Space-time finite elements are particularly advantageous for the
use of concepts such as duality and goal-oriented adaptivity in space and
time~\cite{schmichAdaptivityDynamicMeshes2008,bauseFlexibleGoalorientedAdaptivity2021,rothTensorProductSpaceTimeGoalOriented2023}.
The concepts of variational space-time discretization also provide a unified
approach to stability and error analysis as shown
in~\cite{matthiesHigherOrderVariational2011}.

While this paper deals with tensor-product space-time finite element methods,
there are other approaches to space-time finite elements, particularly for
perturbed space-time meshes. These are discussed in detail
in~\cite{langerSpaceTimeMethodsApplications2019}
and~\cite{steinbachSpaceTimeFiniteElement2015,drflerSpaceTimeDiscontinuousGalerkin2016,langerSpaceTimeHexahedralFinite2022,nochettoSpacetimeMethodsTimedependent2018,ernestiSpaceTimeDiscontinuousPetrov2019a,langerSpaceTimeMethodsApplications2019,langerSpaceTimeHexahedralFinite2022}.
These papers highlight the flexibility and effectiveness of the method in
dealing with complex, time-dependent problems. Other notable examples include
the work
of~\cite{banjaiTrefftzPolynomialSpaceTime2017,gopalakrishnanMappedTentPitching2017},
as well as more recent developments such as those presented
in~\cite{perugiaTentPitchingTrefftzDG2020,steinbachCoerciveSpacetimeFinite2020}.
These works and their references serve as a good basis for a comprehensive
survey of recent developments in space-time discretization techniques.

While time parallelism arises naturally in the context of space-time finite
element methods, parallel time integration is a much broader topic. Most methods
are based on time-domain decomposition techniques and rely on
predictor-corrector algorithms, multilevel, multigrid, or multi-shooting
methods. See~\cite{gander50YearsTime2015,ongApplicationsTimeParallelization2020}
for surveys on time parallelism. In particular, we mention work related to
multigrid
methods~\cite{falgoutParallelTimeIntegration2014,dnnebackeIncreasedSpaceparallelismTimesimultaneous2021,desterckConvergenceAnalysisParallelintime2020}.
In addition to time-domain decomposition approaches, there is also stage
parallelism, i.\,e.\, parallelism within a single time step, see for
example~\cite{christliebParallelHighOrderIntegrators2010,paznerStageparallelFullyImplicit2017,munchStageParallelFullyImplicit2023}.
The scalability of these methods is constrained by the number of stages, which
may result in less parallelization potential than domain decomposition
approaches, for example. Despite this limitation, these methods have been shown
to be effective in the scaling limit~\cite{munchStageParallelFullyImplicit2023},
which is the main motivation for time parallelism in the first place.

There are two approaches to extending classical multigrid methods to the time
dimension: multigrid-in-time and space-time multigrid. In recent years,
significant work has been done on the multigrid-reduction-in-time (MGRIT)
algorithm~\cite{falgoutParallelTimeIntegration2014,hahneAsynchronousTruncatedMultigridReductioninTime2023,ganderMultigridInterpretationsParareal2018,tielenCombiningPmultigridMultigrid2022},
which is based on multigrid in time. An advantage of MGRIT is its easy
integration into existing codes. It only requires a routine to integrate from
one time to the next with an adjustable time step. STMG methods, on the other
hand, have the full advantage of computing multiple time steps at
once~\cite{hortonSpaceTimeMultigridMethod1995,francoMultigridMethodBased2018,falgoutMultigridMethodsSpace2017b,borzExperiencesSpaceTime2005,donatelliAllatonceMultigridApproaches2021,honBlockToeplitzPreconditioner2023}.
Note that sequential time stepping is already optimal in terms of the
complexity, but doesn't provide any parallel scalability. Therefore, paying with
computational complexity to get better scalability is a viable option. In
STMG methods, time is simply another dimension in the grid.
Therefore, given a suitable smoother, they are of optimal complexity. The
systems arising from space-time finite element discretizations are well-suited
for solution by STMG methods. However, most work is based on
algebraic multigrid
methods~\cite{steinbachAlgebraicMultigridMethod2018,langerSpaceTimeFiniteElement2020,langerSpaceTimeHexahedralFinite2022}.
Geometric STMG methods have first been addressed
in~\cite{hackbuschParabolicMultigridMethods1985} and later
in~\cite{ganderAnalysisNewSpaceTime2016,hortonSpaceTimeMultigridMethod1995} as
well as in the references therein. In terms of the discretization, this work has
some similarities to~\cite{ganderAnalysisNewSpaceTime2016}, but we take a
different approach. We consider discontinuous and continuous discretizations in
time for heat and wave equations, and we use a different smoother and focus on
more practical aspects of space-time multigrid. Other works use geometric
multigrid methods, but only consider the potential of parallelism within one
time
step~\cite{anselmannGeometricMultigridMethod2023,anselmannEnergyefficientGMRESMultigrid2024}.
The authors
of~\cite{hoferParallelRobustPreconditioning2019} introduce robust
preconditioners for space-time isogeometric analysis of parabolic evolution
problems, utilizing a time-parallel multigrid approach with discontinuous
variational time discretizations, which we also employ here. The smoother is
specifically designed for symmetric saddle point problems and exploits the
structure of the operator.

Recently, there has been some contributions towards the all-at-once solution of
PDEs leveraging block Toeplitz structures to enhance computational efficiency
and
parallelizability~\cite{mcdonaldPreconditioningIterativeSolution2018,honBlockToeplitzPreconditioner2023,goddardNoteParallelPreconditioning2019,sunParallelinTimeImplementationNumerov2022,danieliAllatonceSolutionLinear2021}.
In most of these contributions, nonsymmetric block Toeplitz structures are
transformed into symmetric systems. The process of symmetrization facilitates
the application of more efficient and robust mathematical tools and algorithms.
The aforementioned works primarily examine finite difference discretizations of
linear parabolic and hyperbolic problems.

Overall, there is a trend towards leveraging space-time approaches in numerical
solution methods for PDEs. This shift is motivated by the need to efficiently
address complex problems. As modern processors have reached their clock speed
limits, the focus of hardware development is on the increase of the number of
processors and the transition to massively parallel computing. To effectively
utilize these advances, scalable, efficient, and flexible computational methods
are essential. Within this field, one of the central problems that needs to be
addressed is the design of performant preconditioners for the arising large
linear systems of equations. We argue that multigrid is one of the most
promising approaches to addressing these challenges. Then, in the context of
STMG methods, the development of efficient smoothers represents another key area
of research. The proposed space-time cell-wise ASM smoother turns
out to be one of the key features in the design of preconditioners for large
linear systems arising from space-time discretizations of PDEs of different
type.

The rest of this paper is structured as follows. Section~\ref{sec:org73458ad}
introduces the mathematical notation and preliminaries essential for the
subsequent discussions. In Section~\ref{sec:orgfee5221} we introduce the
space-time finite element discretization, with an emphasis on the continuous and
discontinuous approaches to time discretization. In Section~\ref{sec:orgbc27d1f}
we introduce the STMG method including the smoother we use to
solve the arising linear systems. The discretization methods are verified
through convergence tests, and the performance of the linear solvers is assessed
via comprehensive scaling experiments in~\ref{sec:org6436fcc}.

\subsection{Notation}
\label{sec:org73458ad}
Let \(\Omega \subset \mathbb{R}^d\) for \(d \in \{1, 2, 3\}\) be a bounded
domain with boundary \(\partial \Omega = \Gamma_D\) and let \(I = (0, T]\)
denote a bounded time interval with final time \(T > 0\). We split the time
interval \(I\) into a sequence of \(N_I\) disjoint subintervals
\(I_n=(t_{n-1},\,t_n]\), \(n=1,\dots,\,N_I\) with the time step
$\tau_n\coloneq t_n-t_{n-1}$. By $\symcal T_{\tau}$ we refer to the time mesh
(or triangulation). For the space discretization, let $\symcal T_h$ be the
quasi-uniform decomposition of $\Omega$ into quadrilaterals or hexahedra with
mesh size $h>0$. Then we denote the space-time tensor product mesh by
$\symcal T_{\tau,\,h}=\symcal T_h\times \symcal T_{\tau}$.

We denote by \(H^1(\Omega)\) the Sobolev space of \(L^2(\Omega)\) functions
whose first derivatives are in \(L^2(\Omega)\). Define
\(H \coloneq L^2(\Omega)\), \(V \coloneq H^1(\Omega)\), and
\(V_0 \coloneq H^1_0(\Omega)\) as the space of \(H^1\)-functions with vanishing
trace on the Dirichlet boundary \(\Gamma_D\). The \(L^2\)-inner product is
\(( \cdot, \cdot )\), with the norm
\(\|\cdot\| \coloneq \|\cdot\|_{L^2(\Omega)}\). %
 For \(J \subseteq [0,\,T]\)
Bochner spaces of \(B\)-valued functions for a Banach space \(B\) are denoted as
\(L^2(J;\,B)\) and \(C(J;\,B)\), each equipped with their natural norms. For a
Banach space \(B\) and \(k\in \mathbb{N}_{0}\) we define
\begin{equation}
  \label{eq:timespace}
  \mathbb{P}_{k}(I_{n};\,B)=\brc{ w_{\tau_{n}}\colon I_{n}\to B\suchthat
    w_{\tau_{n}}(t)=\sum_{j=0}^kW^jt^j\;\forall t\in I_n,\:W^j\in B \;\forall j}\,.
\end{equation}
For \(p\in \mathbb{N}\) we define the finite element space that is built on the
spatial mesh as
\begin{equation}
  \label{eq:fespace}
  \symcal{V}_{h}=\brc{v_h\in C(\bar{\Omega})^d\suchthat v_h\restrict{K} \in
    \mathbb{Q}_p(K)^d\;\forall K \in \mathcal{T}_h}\,, \qquad\symcal{V}_{h,\,0}=\symcal{V}_{h}\cap V_0\,,
\end{equation}
where \(\mathbb{Q}_p(K)\) is the space defined by the reference mapping of
polynomials on the reference element with maximum degree \(p\) in each variable.
For an integer \(k\in \N\) we introduce the space of continuous in time
functions
\begin{equation}
  \label{eq:st-disc-continuous}
  \symcal{X}_{\tau}^k(B)=\brc{w\in C(\bar{I};\,B) \suchthat w\restrict{I_n} \in \mathbb{P}_{k}(I_{n};\,B) \;\forall n=1,\dots,\,N}\,.
\end{equation}
For an integer \(l\in \N_0\) we introduce the space of \(L^2\) in time functions
\begin{equation}
  \label{eq:st-disc-l2}
  \symcal{Y}_{\tau,\,h}^l(B)=\brc{w\in L^2(I;\,B) \suchthat w\restrict{I_n} \in \mathbb{P}_{l}(I_{n};\,B) \;\forall n=1,\dots,\,N}\,.
\end{equation}
For \(B=\symcal{V}_h\) we abbreviate the discrete space-time function spaces as
follows,
\begin{equation}
  \label{eq:disc-spaces}\symcal{X}_{\tau,\,h}^k=\symcal{X}_{\tau}^k(\symcal{V}_h),\qquad\symcal{Y}_{\tau,\,h}^l=\symcal{Y}_{\tau,\,h}^l(\symcal{V}_h)\,.
\end{equation}
We express the \(w_{\tau,\,h}\restrict{{I}_n} \in \P_k{(I_n;\,\symcal{V}_h)}\)
in terms of space- and time-basis functions \(\phi_j(x)\) and \(\xi_{i}(t)\)
\begin{equation}
  \label{eq:spacetime_discrete}
  w_{\tau,\,h}\restrict{{I}_n}(\symbfit x,\,t) \coloneq \sum_{i=1}^{k+1} w_{n}^i(\symbfit x)\xi_{n,\,i}{(t)}=
  \sum_{i=1}^{k+1} \sum_{j=0}^{N_{\symbfit x}} w_{n,\,j}^i\phi_j(\symbfit x)\xi_{n,\,i}{(t)}\,,\qquad \text{for }(\symbfit x,\,t)\in \Omega\times \bar{I_{n}}.
\end{equation}
where \(w_{n,\,i}\) represents the coefficient function to the \(i\)-th time
basis function on \(I_n\) and \(w_{n,\,j}^i\) coefficients to the space-time
basis functions. This representation reveals the important tensor-product
structure, which is useful for the numerical analysis and the implementation.

\section{Space-time finite element discretization}
\label{sec:orgfee5221}
We consider the heat and wave equation as a prototypical parabolic and
hyperbolic problem. The heat equation with thermal diffusivity $\rho\in L^{2}(\R)$ is given
by
\begin{equation}
  \label{eq:heat}
  \partial_{t} u - \nabla\cdot (\rho\nabla u)= f\,,
\end{equation}
equipped with appropriate initial and boundary conditions. Integrating the
strong formulation in space and time and multiplying by space-time test
functions yields the global space-time variational form: Find \(u \in
W(I)\coloneq\{u\in L^2(I;\,V)\suchthat \partial_t  u\in L^2(I;\,H)\}\) such that
for all \(w \in W_0(I)\coloneq\{ w\in L^2(I;\, V_0)\suchthat \partial_t
w\in L^2(I;\, L)\}\),
\[
\int_I \left( ( \partial_t u,\,w ) + ( \rho \nabla u,\,\nabla w
  ) \right) \,\,\mathrm{d}t = \int_I ( f,\,w ) \drv t + ( u(0),\,w(0) ).
\]

The acoustic wave
equation with sound speed $\rho\in L^{2}(\R)$, written as a first order in time system, is
given by
\begin{equation}
  \label{eq:wave}
  \partial_{t} u -  v= 0,\quad
  \partial_{t} v - \nabla\cdot (\rho\nabla u)= f\,,
\end{equation}
equipped with appropriate initial and boundary conditions. Integrating the
strong formulation in space and time and multiplying by space-time test
functions yields the global space-time variational form: Find \((v, u) \in
\left( W(I) \times W(I)\right)\) such that for all \((\tilde w, w) \in W_0(I) \times W_0(I)\),
\[
\int_I \left( ( \partial_t u,\,\tilde w ) - ( v,\,\tilde w
  ) \right) \drv t = 0,\quad\int_I \left( ( \partial_t
  v,\,w ) + ( \rho \nabla u,\,\nabla w )
\right) \drv t = \int_I ( f,\,w ) \drv t + ( v(0), w(0) ).
\]

For the discretization in space, let
\(\{\phi_j\}_{j=1}^{N_{\symbfit x}}\subset \symcal V_h\) denote a nodal
Lagrangian basis of \(\symcal V_h\). The mass and stiffness matrix \(\symbf M_h\)
and \(\symbf A_h\) are defined by
\begin{equation}
  \label{eq:space-matrices}
  \symbf M_h\coloneq ((\phi_i,\,\phi_j))_{i,\,j=1}^{N_{\symbfit x}},
  \quad \symbf A_h \coloneq  (( \rho\nabla \phi_i,\,\nabla \phi_j))_{i,\,j=1}^{N_{\symbfit x}}\,.
\end{equation}
Two types of time discretizations are considered: $DG(k)$, a discontinuous
Galerkin method, and $CGP(k)$, a continuous Galerkin-Petrov discretization, both
of order $k$. The A-stability of the $CGP(k)$ method for $k \ge 1$ and the
L-stability of the $DG(k)$ method for $k \ge 0$ were established
in~\cite{schieweckStableDiscontinuousGalerkinPetrov}
and~\cite{hairerSolvingOrdinaryDifferential1996}. We describe these
discretization schemes in the following sections.

\subsection{The discontinuous Galerkin time discretization}
\label{sec:orge1064c0}
We introduce the $DG(k)$, \(k \ge 0\) methods for the space-time
discretization, utilizing \(\symcal{Y}_{\tau,\,h}^k\)
(cf.~\eqref{eq:disc-spaces}) for the trial and test spaces. This allows for
discontinuities in the trial functions across each subinterval \(I_n\). For
brevity, we assume that the time steps are of uniform size, i.\,e.\ $\tau=\tau_n$
for $n=1,\dots,\,N_I$. The resulting systems decouple at the endpoints of the
subintervals \(I_n\), facilitating a piecewise solution approach across each
interval. For discontinuous functions \(w_\tau(t) \in \symcal{Y}_{\tau,\,h}^k\),
we define:
\begin{equation}\label{eq:jump-terms}
  w_n^- \coloneq \lim_{t \to t_n^-} w_{\tau,\,h}(t),\quad w_n^+ \coloneq \lim_{t \to t_n^+} w_{\tau,\,h}(t),\quad [w_{\tau,\,h}]_n \coloneq w_n^+ - w_n^-.
\end{equation}
To obtain a time marching problem, we introduce a temporal test basis that is
supported on the subintervals $I_n$. Let
\(\{\xi_i\}_{i=1}^{k+1}\subset \mathbb{P}_k(I_n,\,\R)\) and
\(\{\hat \xi_i\}_{i=1}^{k+1}\subset \mathbb{P}_k(\hat I,\,\R)\),
\(\hat I \coloneq [0,\,1]\) denote the Lagrangian basis of
\(\mathbb{P}_k(I_n,\,\R)\) and \(\mathbb{P}_k(\hat I,\,\R)\) w.\,r.\,t.\ the
integration points of the \(k+1\) point Gauss-Radau quadrature. With the affine
transformation
\begin{equation}
  \label{eq:affine-trans}
  \mat{T}_{n}\colon\hat I \to I_{n},\quad \hat t \mapsto t_{n-1} + (t_n-t_{n-1}) \hat t\,,
\end{equation}
the \(i\)th basis function \(\xi_{i}\) on \(I_{n}\) is given by the composition
of
\begin{equation}
  \label{eq:affine-trans-composition}
  \hat\xi_{i}\circ \mat{T}_{n}^{-1}\eqcolon \xi_{i}\,,
\end{equation}
for \(i=1,\dots,\,k+1\). Then, the matrices for the time discretization are
given through the weights \(\symbf M_h\in \R^{N_t\times N_t}\),
\(\symbf A_\tau\in \R^{N_t\times N_t}\) and \(\symbf\alpha\in \R^{N_t}\), with
\begin{equation}
  \label{eq:dg-time-weights}
  {(\symbf M_\tau)}_{i,\,j} \coloneq
  \tau\int_{\hat I}{\hat\xi_j{(\hat t)}
    \hat \xi_i{(\hat t)} \drv \hat t}
  \,,\quad
  {(\symbf A_\tau)}_{i,\,j} \coloneq \int_{\hat I} {\hat\xi_j^\prime{(\hat t)}
    \hat \xi_i{(\hat t)} \drv \hat t}+ \hat\xi_j^\prime{(1)} \hat \xi_i{(1)}
  \,,\quad
  \symbf\alpha_{i} \coloneq  \hat\xi_j^\prime{(0)} \hat \xi_i{(0)}
  \,,\quad i,\,j=1,\dots,\,N_t\,,
\end{equation}
where $N_t=k+1$ is the number of temporal degrees of freedom. In the following
we directly derive the local subproblems of the time marching problem and then
derive subproblems containing multiple time steps. For the derivation of the
global space-time variational formulations and the resulting global systems of
equations we refer to~\cite{kcherVariationalSpaceTime2014}.
\subsubsection{Discontinuous Galerkin time discretization of the heat equation}
\label{sec:orgbcf6bb5}
Consider the local problem on the interval \(I_n\) where the trajectories
\(u_{\tau,\,h}(t)\) have already been computed for all \(t \in [0,\,t_{n-1}]\)
with initial conditions \(u_{\tau,\,h}(0)=u_{0,\,h}\), where
\({u}_{0,\,h} \coloneq I_{\symcal V_h} u_0\) denotes the projection of
\(v_0 \in V\) onto \(\symcal V_h\). Given
\(u_{\tau,\,h}(t_{n-1})\in \symcal{V}_h\), the variational problem for the
DG(\(k\)) method is formulated as finding
\(u_{\tau,\,h}\in\mathbb{P}_{k}(I_{n};\,\symcal{V}_{h})\), such that for all
\(w_{\tau,\,h} \in \mathbb{P}_{k}(I_{n};\,\symcal{V}_{h})\)
\begin{align}\label{eq:var-heat-dg}
  \int_{I_n} \left( \partial_t u_{\tau,\,h},\, w_{\tau,\,h} \right) + \left( \nabla u_{\tau,\,h},\, \nabla w_{\tau,\,h} \right) \drv t + \left( [u_{\tau,\,h}]_{n-1}, {w}_{n-1}^+ \right) = \int_{I_n} \left( f_{\tau,\,h}, w_{\tau,\,h} \right) \drv t\,,
\end{align}
with $u_{n-1}^{-}=u_{\tau,\,h}\restrict{I_{n-1}}(t_{n-1})$. As the trial and
test functions \( w_{\tau, h} \in \mathcal{Y}_{\tau, h}^k \) can be
discontinuous at interval boundaries, the right and left limits, denoted by
\( w_n^+ \) and \( w_n^- \), as defined in~\eqref{eq:jump-terms}, are not equal.
This gives rise to the jump terms introduced in~\eqref{eq:jump-terms} within the
variational formulation~\eqref{eq:var-heat-dg}. They are essential for
maintaining the conservation properties and the correct flux balances at the
interval boundaries. Upon performing numerical integration, the arising
algebraic system in the $n$-th time step of the \(CG(p)-DG(k)\) space-time
finite element discretization amounts to determining
\(\vct u_n=(\vct u_n^1,\dots,\,\vct u_n^{N_t})^{\top}\),
\(\vct u_n^\bullet \in \R^{{N_{\symbfit x}}}\) such that
\begin{equation}
  \label{eq:heat-cg-dg}
  \underbrace{(\symbf{M}_\tau\otimes \symbf{A}_h + \symbf{A}_\tau \otimes
    \symbf{M}_h)}_{\coloneq \symbf S} \vct u_n= \symbf{M}_\tau\otimes \symbf{M}_h \vct f_n+\symbf{\alpha}\otimes \symbf{M}_h \vct u_{n-1}^{N_t}\,,
\end{equation}
where \(\vct f_n=(\vct f_n^1,\dots,\,\vct f_n^{N_t})^{\top}\),
\(\vct f_{n,\,j}^i=( f(\symbfit x,\,t_n^q),\,w_{j})\). Here, \(t_n^q\),
\(q=1,\dots,k+1\) is the \(q\)-th point of the Gauss-Radau quadrature rule with
\(k+1\) points. In~\eqref{eq:heat-cg-dg}, the blocks of $\vct u_n$ correspond to
the temporal degrees of freedom and the elements of the blocks,
\(\vct u_n^\bullet \in \R^{{N_{\symbfit x}}}\) correspond to the spatial degrees
of freedom. Let \(\symbf 1_l\in \R^{m\times n}\) be the matrix with entries
\(a_{i,\,j}=\delta_{l,j}\). Let
\(\symbf{B}\coloneq \symbf 1_{k+1}\otimes\symbf \alpha\otimes\symbf{M}_h\).
Later we want to apply \(h\)-multigrid in time (we may refer to it as
\(\tau\)-multigrid), so we collect consecutive time steps $n_1,\dots,\,n_{c}$ in
one system for a \(CG(p)-DG(k)\) space-time finite element discretization.
Motivated by~\cite{ganderAnalysisNewSpaceTime2016}, similar formulations has
also been used in~\cite{fischerMOReDWRSpacetime2024}.
Using~\eqref{eq:heat-cg-dg} we derive the algebraic system for the solution of
multiple time steps at once as
\begin{equation}
  \label{eq:heat-cg-dg-at-once}
  \begin{pmatrix}
    \symbf{S}                 &                                         &&\\
    -\symbf{B} & \symbf{S}                 &&\\
                              & \ddots & \ddots &\\
                              &&-\symbf{B} & \symbf{S} &\\
                              &&&-\symbf{B} & \symbf{S}
  \end{pmatrix}\begin{pmatrix}
                 \vct u_{n_1}\\
                 \vct u_{n_2}\\
                 \vdots\\
                 \vct u_{n_{c}-1}\\
                 \vct u_{n_{c}}
               \end{pmatrix}=
               \begin{pmatrix}
                 \symbf{M}_\tau\otimes \symbf{M}_h \vct f_{n_1}+\symbf{\alpha}\otimes \symbf{M}_h \vct  u_{n_1-1}^{N_t}\\
                 \symbf{M}_\tau\otimes \symbf{M}_h \vct f_{n_2}\\
                 \vdots\\
                 \symbf{M}_\tau\otimes \symbf{M}_h \vct f_{n_{c}-1}\\
                 \symbf{M}_\tau\otimes \symbf{M}_h \vct f_{n_{c}}
               \end{pmatrix}\,.
\end{equation}

\subsubsection{Discontinuous Galerkin time discretization of the wave equation}
\label{sec:org77adc7e}
Consider the local problem on the interval \(I_n\) where the trajectories
\(u_{\tau,\,h}(t)\) and \(v_{\tau,\,h}(t)\) have already been computed for all
\(t \in [0,\,t_{n-1}]\) with initial conditions \(u_{\tau,\,h}(0)=u_{0,\,h}\),
\(v_{\tau,\,h}(0)=v_{0,\,h}\), where
\({u}_{0,\,h} \coloneq I_{\symcal V_h} u_0\),\,
\({v}_{0,\,h} \coloneq I_{\symcal V_h} v_0\). Given \(u_{\tau,\,h}(t_{n-1})\),
\(v_{\tau,\,h}(t_{n-1})\in \symcal{V}_h\), the variational problem for the
DG(\(k\)) method is formulated as finding
\(u_{\tau,\,h},\,v_{\tau,\,h} \in\mathbb{P}_{k}(I_{n};\,\symcal{V}_{h})\), such
that
\begin{equation}\label{eq:var-wave-dg}
  \begin{aligned}
    \int_{I_n} \left( \partial_t u_{\tau,\,h},\, {\tilde w}_{\tau,\,h} \right) \drv t + \left( [u_{\tau,\,h}]_{n-1},\, {\tilde w}_{n-1}^+ \right)
    - \int_{I_n} \left( v_{\tau,\,h},\, {\tilde w}_{\tau,\,h} \right) \drv t &= 0\,,\\
    \int_{I_n} \left( \partial_t v_{\tau,\,h},\, {w}_{\tau,\,h} \right) \drv t + \left( [v_{\tau,\,h}]_{n-1},\, {w}_{n-1}^+ \right)
    + \int_{I_n} \left( \rho\nabla u_{\tau,\,h},\, \nabla w_{\tau,\,h} \right)\drv t &= \int_{I_n} \left( f,\, {w}_{\tau,\,h} \right) \drv t\,,
  \end{aligned}
\end{equation}
for all
\({w}_{\tau,\,h},\, {\tilde w}_{\tau,\,h} \in
\mathbb{P}_{k}(I_{n};\,\symcal{V}_{h})\), with
$u_{n-1}^{-}=u_{\tau,\,h}\restrict{I_{n-1}}(t_{n-1})$,
$v_{n-1}^{-}=v_{\tau,\,h}\restrict{I_{n-1}}(t_{n-1})$. Note the jump terms
arising in~\eqref{eq:var-wave-dg}, analogous to~\eqref{eq:var-heat-dg}.
Specifically, the jump terms arise for those variables on which a time
derivative is present, see~\cite{kcherVariationalSpaceTime2014} for details.
Upon numerical integration, the resulting algebraic system in the
$n$-th time step of the \(CG(p)-DG(k)\) space-time finite element discretization
amounts to determining \(\vct u_n=(\vct u_n^1,\dots,\,\vct u_n^{N_t})^{\top}\),
\(\vct v_n=(\vct v_n^1,\dots,\,\vct v_n^{N_t})^{\top}\), where
\(\vct u_n^\bullet,\,\vct v_n^\bullet \in \R^{{N_{\symbfit x}}}\), such that
\begin{subequations}
  \label{eq:wave-cg-dg}
  \begin{align}\label{eq:wave-cg-dg-ade}
    \vct v_n&=\symbf{M}_\tau^{-1}\symbf{A}_\tau\vct u_n - \symbf{M}_\tau^{-1}\symbf{\alpha}\vct u_{n-1}^{N_t}\,,\\
    \label{eq:wave-cg-dg-pde}
    \underbrace{(\symbf{M}_\tau\otimes \symbf{A}_h +
    \symbf{A}_\tau \symbf{M}_\tau^{-1}\symbf{A}_\tau \otimes \symbf{M}_h)}_{\coloneq
    \symbf S}\vct u_n&= \symbf{M}_\tau\otimes \symbf{M}_h \vct
                       f+\symbf{\alpha}\otimes \symbf{M}_h\vct v_{n-1}^{N_t} +
                       \underbrace{\symbf{A}_\tau \symbf{M}_\tau^{-1}\symbf\alpha\otimes
                       \symbf{M}_h}_{\coloneq\symbf{D}}\symbf u_{n-1}^{N_t}\,.
  \end{align}
\end{subequations}
In~\eqref{eq:wave-cg-dg-ade}, we observe that the system of equations comprises
solely time derivatives. This allows for the condensation of the system, such
that we only solve for the displacement $\vct u_n$. The velocity $\vct v_n$ is
computed by means of simple vector update equations in a second step. To be able
to use \(\tau\)-multigrid, we collect consecutive time steps $n_1,\dots,\,n_{c}$
in one system for a \(CG(p)-DG(k)\) space-time finite element discretization. To
this end, let
\(\symbf{B}\coloneq{(\symbf M_\tau)}_{k+1,\,k+1}^{-1}\symbf\alpha\symbf {(A_\tau)}_{k+1,\,\cdot}\otimes
\symbf M_h+\symbf{D}\),
\(\symbf{C}\coloneq \symbf {(M_\tau)}_{k+1,\,k+1}^{-1}\symbf
1_{k+1}\otimes\symbf\alpha\otimes \symbf M_h\). Using~\eqref{eq:wave-cg-dg} we
derive the algebraic system for the solution of multiple time steps at once as
\begin{equation}
  \label{eq:wave-cg-dg-at-once}
  \begin{pmatrix}
    \symbf{S}                 &                                         &&\\
    -\symbf{B} & \symbf{S}                 &&\\
    -\symbf{C} &-\symbf{B} & \symbf{S} &\\
    \;\ddots & \;\ddots& \;\ddots &\\
                              &-\symbf{C}&-\symbf{B} & \symbf{S} &\\
                              &&-\symbf{C}&-\symbf{B} & \symbf{S}
  \end{pmatrix}\begin{pmatrix}
                 \vct u_{n_1}\\
                 \vct u_{n_2}\\
                 \vct u_{n_3}\\
                 \vdots\\
                 \vct u_{n_{c}-1}\\
                 \vct u_{n_{c}}
               \end{pmatrix}=
               \begin{pmatrix}
                 \symbf{M}_\tau\otimes \symbf{M}_h \vct f_{n_1}+\symbf{\alpha}\otimes \symbf{M}_h \vct v_{n_1-1}^{N_t}+\symbf{D} \vct u_{n_1-1}^{N_t}\\
                 \symbf{M}_\tau\otimes \symbf{M}_h \vct f_{n_2}+\symbf{C}\vct u_{n_1-1}^{N_t}\\
                 \symbf{M}_\tau\otimes \symbf{M}_h \vct f_{n_3}\\
                 \vdots\\
                 \symbf{M}_\tau\otimes \symbf{M}_h \vct f_{n_{c}-1}\\
                 \symbf{M}_\tau\otimes \symbf{M}_h \vct f_{n_{c}}
               \end{pmatrix}\,.
             \end{equation}

\subsection{The continuous Galerkin-Petrov time discretization}
\label{sec:disc-cgp}
Given the continuous time-discrete trial space \(\symcal{X}_{\tau,\,h}^{k}\)
(cf.~\eqref{eq:disc-spaces}), and the discontinuous time-discrete test space
\(\symcal{Y}_{\tau,\,h}^{k-1}\), the method solves local problems on each
subinterval \(I_n\) with given initial conditions. This is due to the test
functions in \(\symcal{Y}_{\tau,\,h}^{k-1}\) having one polynomial degree less
than the trial functions, which introduces a system of equations that decouples
at \(t_n\) due to the discontinuity of the test functions. The additional degree
of freedom in the trial functions is fixed by the continuity of
$\symcal{X}_{\tau,\,h}^{k}$: We impose a continuity constraint at the left
boundary of \(I_n\), enabling the problem to be solved using a time-marching
process. To this end we introduce a temporal test basis that is supported on the
subintervals $I_n$. Let \(\{\xi_i\}_{i=1}^{k+1}\subset \mathbb{P}_k(I_n,\,\R)\),
\(\{\hat \xi_i\}_{i=1}^{k+1}\subset \mathbb{P}_k(\hat I,\,\R)\),
\(\hat I \coloneq [0,\,1]\) denote the Lagrangian basis of
\(\mathbb{P}_k(I_n,\,\R)\), \(\mathbb{P}_k(\hat I,\,\R)\) w.\,r.\,t.\ the
integration points of the \(k+1\) point Gauss-Lobatto quadrature. For the trial
space, let \(\{\psi_i\}_{i=1}^{k}\subset \mathbb{P}_{k-1}(I_n,\,\R)\),
\(\{\hat \psi_i\}_{i=1}^{k}\subset \mathbb{P}_{k-1}(\hat I,\,\R)\) denote the
Lagrangian basis of \(\mathbb{P}_{k-1}(I_n,\,\R)\),
\(\mathbb{P}_{k-1}(\hat I,\,\R)\) w.\,r.\,t.\ the last $k$ integration points of
the \(k+1\) point Gauss-Lobatto quadrature.

With the affine transformation~\eqref{eq:affine-trans}, the \(i\)th basis
function \(\xi_{i}\) or \(\psi_i\) on \(I_{n}\) is given
by~\eqref{eq:affine-trans-composition}. Then, the matrices for the time
discretization are given through the weights
\(\symbf M_\tau\in \R^{N_t\times N_t}\), \(\symbf A_\tau\in \R^{N_t\times N_t}\),
\(\symbf\beta\in \R^{N_t}\) and \(\symbf\alpha\in \R^{N_t}\), with
\begin{equation}
  \label{eq:cg-time-weights}
  \begin{aligned}
    {(\symbf M_\tau)}_{i,\,j-1} &\coloneq
                             \tau\int_{\hat I}{\hat\xi_j{(\hat t)}
                             \hat \psi_i{(\hat t)} \drv \hat t}
                             \,,\quad
                             {(\symbf A_\tau)}_{i,\,j-1} \coloneq \int_{\hat I} {\hat\xi_j^\prime{(\hat t)}
                             \hat \psi_i{(\hat t)} \drv \hat t}\,,\\
    \symbf\beta_{i} &\coloneq
                       \tau\int_{\hat I}{\hat\xi_1{(\hat t)}
                       \hat \psi_i{(\hat t)} \drv \hat t}
                       \,,\quad
                       \symbf\alpha_{i} \coloneq \int_{\hat I} {\hat\xi_1^\prime{(\hat t)}
                       \hat \psi_i{(\hat t)} \drv \hat t}\,,\quad i=1,\dots,\,N_t,\; j=2,\dots,\,N_t+1
  \end{aligned}
\end{equation}
where $N_t=k$ is the number of temporal degrees of freedom in the algebraic
system we derive next, as we already separated the weights associated to the degree
of freedom fixed by the continuity constraints into $\symbf\beta$ and
$\symbf\alpha$.

\subsubsection{Continuous Galerkin-Petrov time discretization of the heat equation}
\label{sec:orgc1f0a88}
Consider the local problem on the interval \(I_n\) where the trajectories
\(u_{\tau,\,h}(t)\) have already been computed for all \(t \in [0,\,t_{n-1}]\)
with initial conditions \(u_{\tau,\,h}(0)=u_{0,\,h}\), where
\({u}_{0,\,h} \coloneq I_{\symcal V_h} u_0\) denotes the interpolation of \(v_0 \in H\) onto
\(\symcal V_h\). Given \(u_{\tau,\,h}(t_{n-1})\in \symcal{V}_h\), the
variational problem for the $CGP(k)$ method is formulated as finding
\(u_{\tau,\,h}\in\mathbb{P}_{k}(I_{n};\,\symcal{V}_{h})\), such that for all
\({w}_{\tau,\,h} \in \mathbb{P}_{k-1}(I_{n};\,\symcal{V}_{h})\)
\begin{equation}
  \int_{I_n} \left( \partial_t u_{\tau,\,h}, w_{\tau,\,h} \right) + \left( \nabla u_{\tau,\,h}, \nabla w_{\tau,\,h} \right) \drv t = \int_{I_n} \left( f, w_{\tau,\,h} \right) \drv t\,,
\end{equation}
with the continuity constraint
$u_{n-1}^{-}=u_{\tau,\,h}\restrict{I_{n-1}}(t_{n-1})$. Upon numerical
integration (cf.~\cite{kcherVariationalSpaceTime2014}), the algebraic system
that arises in the $n$-th time step of the \(CG(p)-CGP(k)\) space-time finite
element discretization amounts to determining
\(\vct u_n=(\vct u_n^1,\dots,\,\vct u_n^{N_t})^{\top}\), where
\(\vct u_n^\bullet \in \R^{{N_{\symbfit x}}}\) such that
\begin{equation}
  \label{eq:heat-cg-cg}
  \underbrace{(\symbf{M}_\tau\otimes \symbf{A}_h +
    \symbf{A}_\tau \otimes \symbf{M}_h)}_{\coloneq \symbf{S}} \vct u_n= \symbf{M}_\tau\otimes \symbf{M}_h \vct f_n - \symbf{\beta}\otimes \symbf{M}_h \vct f_{n-1}^{N_t} + \underbrace{(\symbf{\beta}\otimes \symbf{A}_h + \symbf{\alpha} \otimes \symbf{M}_h)}_{\coloneq \symbf{b}}\vct u_{n-1}^{N_t}
\end{equation}
where \(\vct f_n=(\vct f_n^1,\dots,\,\vct f_n^{N_t})^{\top}\),
\(\vct f_{n,\,j}^i=( f(\symbfit x,\,t_n^q),\,w_{j})\). The points
\(t_n^q\), \(q=1,\dots,k\) are the \(q+1\)-th nodes of the Gauss-Lobatto
quadrature rule with \(k+1\) points. Note that by this choice of quadrature rule
the continuity constraint can be written as $\vct u_{n-1}^{N_t}=\vct u_{n}^{0}$.
Here, $\vct u_{n}^{0}$ is the temporal degree of freedom at the first
Gauss-Lobatto quadrature point, with the zero index indicating that it acts
algebraically like an initial condition. Let
\(\symbf{B}\coloneq\symbf 1_k\otimes\symbf{b}\). Analogous to the DG time stepping, we collect
consecutive time steps \(n_1,\dots,\,n_{c}\) in one system for a \(CG(p)-CGP(k)\)
space-time finite element discretization. Using~\eqref{eq:heat-cg-cg} we derive
the algebraic system for the solution of multiple time steps at once as
\begin{equation}
  \label{eq:heat-cg-cg-at-once}
  \begin{pmatrix}
    \symbf{S}                 &                                         &&\\
    -\symbf{B} & \symbf{S}                 &&\\
                              & \ddots & \ddots &\\
                              &&-\symbf{B} & \symbf{S} &\\
                              &&&-\symbf{B} & \symbf{S}
  \end{pmatrix}
  \begin{pmatrix}
    \vct u_{n_1}\\
    \vct u_{n_2}\\
    \vdots\\
    \vct u_{n_{c}-1}\\
    \vct u_{n_{c}}
  \end{pmatrix}=
  \begin{pmatrix}
    \symbf{M}_\tau\otimes \symbf{M}_h \vct f_{n_1} - \symbf{\beta}\otimes \symbf{M}_h \vct f_{n_1-1}^{N_t} + \symbf{b} u_{n_1-1}^{N_t}\\
    \symbf{M}_\tau\otimes \symbf{M}_h\vct f_{n_2} - \symbf{\beta}\otimes \symbf{M}_h\vct f_{n_{1}}^{N_t}\\
    \vdots\\
    \symbf{M}_\tau\otimes \symbf{M}_h\vct f_{n_{c}-1} - \symbf{\beta}\otimes \symbf{M}_h\vct f_{n_{c}-2}^{N_t}\\
    \symbf{M}_\tau\otimes \symbf{M}_h\vct f_{n_{c}} - \symbf{\beta}\otimes \symbf{M}_h\vct f_{n_{c}-1}^{N_t}
  \end{pmatrix}\,.
\end{equation}

\subsubsection{Continuous Galerkin-Petrov time discretization of the wave equation}
\label{sec:org05b7863}
Consider the local problem on the interval \(I_n\) where the trajectories
\(u_{\tau,\,h}(t)\) and \(v_{\tau,\,h}(t)\) have already been computed for all
\(t \in [0,\,t_{n-1}]\) with initial conditions \(u_{\tau,\,h}(0)=u_{0,\,h}\),
\(v_{\tau,\,h}(0)=v_{0,\,h}\), where
\({u}_{0,\,h} \coloneq I_{\symcal V_h} u_0\),\,
\({v}_{0,\,h} \coloneq I_{\symcal V_h} v_0\). Given \(u_{\tau,\,h}(t_{n-1})\),
\(v_{\tau,\,h}(t_{n-1})\in \symcal{V}_h\), the variational problem for the
DG(\(k\)) method is formulated as finding
\(u_{\tau,\,h},\,v_{\tau,\,h} \in\mathbb{P}_{k-1}(I_{n};\,\symcal{V}_{h})\),
such that
\begin{equation}
  \int_{I_n}{\left(
      \partial_t u_{\tau,\,h}, {w}_{\tau,\,h}
    \right)
    - \left( v_{\tau,\,h}, {w}_{\tau,\,h} \right)
    \,\drv t } = 0,\quad
  \int_{I_n}{\left(
      \partial_t v_{\tau,\,h}, {w}_{\tau,\,h}
    \right)
    + \left( \rho\nabla u_{\tau,\,h},\, \nabla w_{\tau,\,h} \right)
    \,\drv t }
  = \int_{I_n}{
    \left( f, {w}_{\tau,\,h} \right) \,\drv t },
\end{equation}
for all
\({w}_{\tau,\,h},\, {\tilde w}_{\tau,\,h} \in
\mathbb{P}_{k}(I_{n};\,\symcal{V}_{h})\), with
$u_{n-1}^{-}=u_{\tau,\,h}\restrict{I_{n-1}}(t_{n-1})$,
$v_{n-1}^{-}=v_{\tau,\,h}\restrict{I_{n-1}}(t_{n-1})$. Upon performing numerical
integration, the arising algebraic system in the $n$-th time step of the
\(CG(p)-CGP(k)\) space-time finite element discretization amounts to determining
\(\vct u_n=(\vct u_n^1,\dots,\,\vct u_n^{N_t})^{\top}\),
\(\vct v_n=(\vct v_n^1,\dots,\,\vct v_n^{N_t})^{\top}\),
\(\vct u_n^\bullet,\,\vct v_n^\bullet \in \R^{{N_{\symbfit x}}}\) such that
\begin{subequations}
\label{eq:wave-cg-cg}
\begin{align}
  \label{eq:wave-cg-cg-ade}
    \vct v_n&=\symbf{M}_\tau^{-1}\symbf{A}_\tau\vct u_n - \symbf{M}_\tau^{-1}\symbf{\alpha}\vct u_{n-1}^{N_t} + \symbf{M}_\tau^{-1}\symbf{\beta}\vct v_{n-1}^{N_t}\,,\\[3pt]
  \begin{split}\label{eq:wave-cg-cg-system}
\underbrace{(\symbf{M}_\tau\otimes \symbf{A}_h + \symbf{A}_\tau \symbf{M}_\tau^{-1}\symbf{A}_\tau \otimes \symbf{M}_h)}_{\coloneq \symbf{S}}\vct u_n&= \symbf{M}_\tau\otimes \symbf{M}_h \vct f - \symbf{\beta}\otimes \symbf{M}_h\vct f_{n-1}^{N_t}\\[-14pt]
 &+ (\symbf{\beta}\otimes \symbf{A}_h+\symbf{A}_\tau \symbf{M}_\tau^{-1}\symbf{\alpha} \otimes \symbf{M}_h)\vct u_{n-1}^{N_t} + (\symbf{\alpha}-\symbf{A}_\tau \symbf{M}_\tau^{-1}\symbf\beta) \otimes \symbf{M}_h\vct v_{n-1}^{N_t}\,.
  \end{split}
\end{align}
\end{subequations}
Analogous to the DG time stepping, we condense the system and compute $\vct v_n$
by vector update equations and only solve for the displacement $\vct u_n$. We
further collect consecutive time steps \(n_1,\dots,\,n_{c}\) in one system for a
\(CG(p)-CGP(k)\) space-time finite element discretization of the wave equation.
The auxiliary variable $\vct v_{n}$ complicates the system because it is defined
recursively by~\eqref{eq:wave-cg-cg-ade} and depends on the past solution. As
for the discontinuous time stepping, this recurrence introduces additional terms
in the matrix blocks. Due to the continuity constraints on both variables, the
system becomes significantly more complex, at least if we don't introduce higher
order derivatives and still want to solve for $\vct v$ through vector updates.
On the right hand side of~\eqref{eq:wave-cg-cg-system} we observe that, in order
to derive the system for the solution of multiple time steps at once, we need an
equation for the initial value $\vct v_{l-1}^{N_t}$ on each interval $I_l$ for
$l\in \{n_2,\dots,\,n_{c}\}$. Note that we can see
in~\eqref{eq:wave-cg-cg-at-once}, for $l=n_1$, we simply
recover~\eqref{eq:wave-cg-cg-system} and $\vct v_{n_1-1}^{N_t}$ is known from
previous computations. From~\eqref{eq:wave-cg-cg-ade} we observe that, in
contrast to the discontinuous time discretization, the initial value
$\vct v_{l-1}^{N_t}$, $l\in \{n_2,\dots,\,n_{c}\}$ cannot be expressed solely in
terms of the displacement variable. For instance,
considering~\eqref{eq:wave-cg-cg-system} for $n=n_2$, we see that we need to
replace $\vct v_{n_1}^{N_t}$ on the right hand side by deriving an equation for
$\vct v_{n_1}^{N_t}$ in terms of the variables $\vct u_{n_1}$ and
$\vct v_{n_1-1}$. Doing this by using the auxiliary differential
equation~\eqref{eq:wave-cg-cg-ade} and iterating this process for all time steps
$n_2,\dots,\,n_{c}$, we get the algebraic system for the solution of multiple
time steps at once as
\begin{subequations}
  \label{eq:wave-cg-cg-at-once}
  \begin{equation}
    \begin{pmatrix}
      \symbf{S}                 &                                         &&\\
      -\symbf{B}+\symbf{E}_{2,\,1} & \symbf{S}                 &&\\
      \symbf{E}_{3,\,1} &-\symbf{B}+\symbf{E}_{3,\,2} & \symbf{S} &\\
                                &\vdots & \;\ddots& \;\ddots &\\
      \vdots&\ddots&-\symbf{B}+\symbf{E}_{n_{c}-1,n_{c}-2} & \symbf{S} &\\
      \symbf{E}_{n_{c},\,1}&\symbf{E}_{n_{c},\,2}&\cdots&-\symbf{B}+\symbf{E}_{n_{c},n_{c}-1} & \symbf{S}
    \end{pmatrix}\begin{pmatrix}
                   \vct u_{n_1}\\
                   \vct u_{n_2}\\
                   \vct u_{n_3}\\
                   \vdots\\
                   \vct u_{n_{c}-1}\\
                   \vct u_{n_{c}}
                 \end{pmatrix}=\mat F
\end{equation}

\begin{equation}
  \mat F =\begin{pmatrix}
            \symbf{M}_\tau\otimes \symbf{M}_h \vct f_{n_1} - \symbf{\beta}\otimes \symbf{M}_h \vct f_{n_1-1}^{N_t}
            + \symbf{\beta}\otimes \symbf{A}_h u_{n_1-1}^{N_t}
            + \symbf A_\tau\symbf M_\tau^{-1}\symbf\alpha\otimes\symbf{M}_h \vct u_{n_1-1}^{N_t}
            +\symbf\alpha_{\symbf{A}_\tau}\otimes\symbf{M}_h\vct v_{n_1-1}^{N_t}  \\
            \symbf{M}_\tau\otimes \symbf{M}_h\vct f_{n_2} - \symbf{\beta}\otimes \symbf{M}_h\vct f_{n_{1}}^{N_t}
            - z_{\symbf M_\tau}\symbf\alpha_{\symbf{A}_\tau}\otimes\symbf{M}_h\vct u_{n_1-1}^{N_t}
            - g_{\symbf M_\tau} \symbf\alpha_{\symbf{A}_\tau}\otimes\symbf{M}_h\vct v_{n_1-1}^{N_t}
            \\
            \symbf{M}_\tau\otimes \symbf{M}_h\vct f_{n_3} - \symbf{\beta}\otimes \symbf{M}_h\vct f_{n_2}^{N_t}
            - z_{\symbf M_\tau} g_{\symbf M_\tau} \symbf\alpha_{\symbf{A}_\tau}\otimes\symbf{M}_h\vct u_{n_1-1}^{N_t}
            - g_{\symbf M_\tau} ^{2}\symbf\alpha_{\symbf{A}_\tau}\otimes\symbf{M}_h\vct v_{n_1-1}^{N_t}
            \\
            \vdots\\
            \symbf{M}_\tau\otimes \symbf{M}_h\vct f_{n_{c}-1} - \symbf{\beta}\otimes \symbf{M}_h\vct f_{n_{c}-2}^{N_t}
            - z_{\symbf M_\tau} g_{\symbf M_\tau} ^{n_{c}-3}\symbf\alpha_{\symbf{A}_\tau}\otimes\symbf{M}_h\vct u_{n_1-1}^{N_t}
            - g_{\symbf M_\tau} ^{n_{c}-2}\symbf\alpha_{\symbf{A}_\tau}\otimes\symbf{M}_h\vct v_{n_1-1}^{N_t}
            \\
            \symbf{M}_\tau\otimes \symbf{M}_h\vct f_{n_{c}} - \symbf{\beta}\otimes \symbf{M}_h\vct f_{n_{c}-1}^{N_t}
            - z_{\symbf M_\tau} g_{\symbf M_\tau} ^{n_{c}-2}\symbf\alpha_{\symbf{A}_\tau}\otimes\symbf{M}_h\vct u_{n_1-1}^{N_t}
            - g_{\symbf M_\tau} ^{n_{c}-1}\symbf\alpha_{\symbf{A}_\tau}\otimes\symbf{M}_h\vct v_{n_1-1}^{N_t}\,,
          \end{pmatrix}
\end{equation}
\end{subequations}
where we used the notation
\(\symbf\alpha_{\symbf{A}_\tau}=\symbf\alpha-\symbf A_\tau\symbf M_\tau^{-1}\symbf\beta\),
\(g_{\symbf M_\tau} =\symbf\beta_{k}{(\symbf M_\tau)}_{k,\,k}^{-1}\) and
\(z_{\symbf M_\tau} =\symbf\alpha_{k}{(\symbf M_\tau)}_{k,\,k}^{-1}\) and
\begin{subequations}
  \begin{align}
    \symbf{B}&\coloneq\symbf
               1_k\otimes\symbf A_\tau\symbf M_\tau^{-1}\symbf\alpha\otimes\symbf{M}_h+\symbf
               1_k\otimes\symbf\beta\otimes\symbf{A}_h\,,\\
    \symbf{E}_{i,\,j}&\coloneq -\symbf{C}_{i,\,j}+\symbf{D}_{i,\,j}\,.
                       \intertext{Here,}
                       \label{eq:c}
                       \symbf{C}_{i,\,j}\coloneq
                       \frac{g_{\symbf M_\tau} ^{i-j-1}}{ {(M_\tau)}_{k,k}}\symbf\alpha_{\symbf{A}_\tau}{(\symbf A_\tau)}_{k,\bullet}\,,&\quad\symbf{D}_{i,\,j}\coloneq
                                                                                                                                                                  \begin{cases}
                                                                                                                                                                    \symbf 1_k\otimes g_{\symbf M_\tau} ^{i-j-2}z_{\symbf M_\tau} \symbf\alpha_{\symbf{A}_\tau}, & i> 2 \and i-1> j \\
                                                                                                                                                                    \symbf 0  &\text{otherwise}
                                                                                                                                                                  \end{cases}\,,
  \end{align}
\end{subequations}
where ${(\symbf A_\tau)}_{k,\bullet}$ is the $k$th row of $\symbf A_\tau$.

\section{Space-time multigrid method}
\label{sec:orgbc27d1f}
We present the space-time multigrid method to solve the space-time
systems~\eqref{eq:heat-cg-dg-at-once},~\eqref{eq:wave-cg-dg-at-once},~\eqref{eq:heat-cg-cg-at-once}
and~\eqref{eq:wave-cg-cg-at-once}. We write the linear systems associated to
these problems as
\begin{equation}
  \label{eq:lin-sys}
  \symbfcal S \vct u=\vct f\,.
\end{equation}
To build the hierarchical sequence of space-time meshes
$\{\symcal T_{\tau,\,h}^l\}_{l=0}^L$, we only coarsen either in space or in
time, as sketched in Figure~\ref{fig:stmg}. Of course the reduction in degrees
of freedom on each level makes space-time coarsening particularly appealing.
Although it is practically feasible, it is challenging in the parallel in time
context, as observed for instance in~\cite{angelImpactSpatialCoarsening2021}.
In~\cite{chaudet-dumasOptimizedSpaceTimeMultigrid2023} the authors derive a
CFL-type condition, which ensures convergence with space-time coarsening for a
one-dimensional STMG method for finite differences
discretizations of the heat equation. Our current investigation is focused on
the computational aspects of the proposed method and the efficacy of its
implementation. Consequently, the exploration of space-time coarsening
techniques is not included in this study and is subject to future work. The
extension within our code is straightforward. Analogous to the notion of $h$
and $p$ multigrid we denote the same things in time by $\tau$ and $k$ multigrid
within this paper. Both modes are implemented in our code.

Each level of the space-time multigrid hierarchy is associated to a system of
linear equations denoted by $\symbfcal S_{l} \vct u_l=\vct f_l$. The finest
level is associated to the original problem~\eqref{eq:lin-sys}. We apply the
multigrid in a $V$-cycle with a space-time cell-based ASM smoother
as preconditioner of GMRES.

\subsection{Multigrid algorithm}
\label{sec:org1919a00}
\begin{figure}[htbp]
  \centering
  \includegraphics[width=\linewidth]{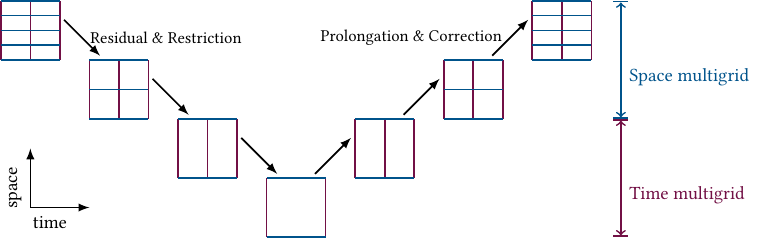}
  \caption{A sketch of the STMG method, see also Algorithm~\ref{alg:stmg}. The
    corrections are transferred by the prolongation operators and the residual
    is transferred by the restriction operators. On each level the error is
    \emph{smoothed} by a single iteration of the additive Schwarz smoother. The
    preferred coarsening strategy, which is then used in
    Algorithm~\ref{alg:mg-sequence}, is first in space and then in
    time.}\label{fig:stmg}
\end{figure}
Space-time multigrid methods expand upon the classical multigrid approach by
incorporating space and time dimensions into the grid coarsening process
(cf.~Figure~\ref{fig:stmg}). The objective of this extension is to efficiently
solve large-scale partial differential equations by performing coarsening and
refinement across space and time dimensions. These methods are of particular
interest when spatial parallelism saturates. To gain further parallelization
potential, the sizes of the subproblems in every time step are increased by
batching multiple time steps in one. Further, when high order discretizations are
employed on problems with large coarse problems, temporal multigrid is useful
for reducing the size of the coarse problem. For an introduction and more
detailed discussion of the geometric multigrid method, we refer
to~\cite{hackbuschMultigridMethodsApplications1985,brambleMultigridMethods1993,trottenbergMultigrid2000,vassilevskiMultilevelBlockFactorization2008}.

\begin{algorithm}[H]
  \caption{Space-time multigrid algorithm
    for problems discretized with space-time finite
    elements~\eqref{eq:lin-sys}. Note the recursive call, which corresponds to a
    $V$-cycle.}\label{alg:stmg}
  \LinesNumbered %
  \SetKwProg{Function}{function}{}{end} %
  \Function{$\operatorname{\mathtt{Multigrid}}(\ell,\;\symbfcal
    S_{\ell},\;\symbf f_{\ell},\symbf u_{\ell})$}{
  \label{alg:presmooth}
  Pre-smooth $\vct u_\ell$ with $\nu_1$ iterations of a
  smoother\;%
  Compute and restrict residual
  $\vct r_{\ell-1} = \symcal R_\ell(\vct f_\ell - \symbfcal S_\ell \symbf u_\ell)$\;
  \uIf{$\ell$ is the coarsest level}{
    Solve $\symbfcal S_\ell \vct u_\ell = \vct f_\ell$\;
  }\Else{
    $\vct e_{\ell-1} = \operatorname{\mathtt{Multigrid}}(l-1,\,\symbfcal
    S_{\ell-1},\,\vct r_{\ell-1},\,\symbf 0$)\;
  }
  Prolongate and correct
  $\vct u_\ell = \vct u_\ell + \symcal P_{\ell} \vct e_{\ell-1}$\;
  \label{alg:postsmooth}
  Post-smooth $\vct u_\ell$ with $\nu_2$ iterations of a
  smoother\;%
  }
\end{algorithm}

The computational complexity of a space-time multigrid cycle
(cf.~Algorithm~\ref{alg:stmg}) is similar to that of multigrid cycles that focus
only on spatial coarsening. In the case of spatial coarsening, the number of
grid points is reduced by a factor of $2^d$ at each level, where $d$ is the
spatial dimension. Similar to a one-dimensional spatial coarsening, the number
of time steps is halved at each coarsening level. For polynomial multigrid the
coarsening factor is derived analogously, depending on the ratio between the
polynomial degrees on each level.

In the previous section, we have already formulated the linear systems, that
arise from the discretization. For the evaluation of operators we utilize
matrix-free methods, where the system matrix is not explicitly formed and stored
and exploit the tensor product structure for the efficient evaluation of finite
element operators by decomposing multidimensional into one-dimensional
operations~\cite{kronbichlerGenericInterfaceParallel2012}.
In matrix-free methods the matrix-vector product \( \vct y = \symbf{S} \vct u \)
is computed directly through local operations and mappings from local to global
degrees of freedom:
\[ \symbf{S} = \sum_{c=1}^{n_c} \symbf{P}_{c,\,\text{loc-glob}}^T \symbf{S}_c
  \symbf{P}_{c,\,\text{loc-glob}}, \] where
\( \symbf{P}_{c,\,\text{loc-glob}} \) maps local degrees of freedom to global
indices, and \( \symbf{S}_c \) is defined by the relation
\( \symbf{S}_c = \symbf{B}_c^T \symbf{D}_c \symbf{B}_c \). If $\symbf{S}$ is a
stiffness matrix, \( \symbf{B}_c \) represents the gradients of shape functions,
and \( \symbf{D}_c \) is a diagonal matrix containing quadrature weights and
coefficients. The integrals are evaluated on-the-fly. The sum factorization
reduces the computational complexity from \(O(p^{2d})\) to \(O(d p^2)\), making
it significantly more efficient for high-dimensional and high-degree polynomial
problems.

The sum-factorization
allows for the fast evaluation by the transformation of multi-dimensional
operations into sequences of one-dimensional operations and reduces the
computational complexity from \( \mathcal{O}({n_{1D}}^{2d}) \) to
\( \mathcal{O}({d n_{1D}^{d+1}}) \), where $n_{1D}$ is the
number of degrees of freedom in one spatial direction. Specifically, the local
matrix-vector product is represented as
\[ \symbf{S}_c \cdot u_c = \symbf{B}_c^T (\symbf{D}_c (\symbf{B}_c \cdot
  u_c)), \] with the evaluation of three matrix-vector products, as indicated by
the braces, which results in the reduced computational complexity. These
techniques are employed for matrix-vector products with the stiffness matrix
$\symbf A_h$ and the mass matrix $\symbf M_h$ in the space-time
systems~\eqref{eq:heat-cg-dg-at-once},~\eqref{eq:wave-cg-dg-at-once},~\eqref{eq:heat-cg-cg-at-once}
and~\eqref{eq:wave-cg-cg-at-once}. The temporal matrices on the other
hand are precomputed and stored as full matrices for each subproblem. As they
are small and $N_{t}\ll N_{\symbfit x}$, this is reasonable and efficient in the
current implementation. Moreover, one can evaluate
$(\symbf{M}_\tau\otimes \symbf{A}_h + \symbf{A}_\tau \otimes \symbf{M}_h)\vct u$
as
$\symbf{M}_\tau (\symbf{A}_h\otimes \vct u) + \symbf{A}_\tau (\symbf{M}_h\otimes
\vct u)$, i.\,e.\ calculating $\symbf{A}_h\vct u^i$ and $\symbf{M}_h\vct u^{i}$
once for each $i=1,\dots,\,N_t$ and reusing it for the multiplication with the
temporal matrices throughout the operator evaluation. This is used analogously
for other Kronecker products
in~\eqref{eq:heat-cg-dg-at-once},~\eqref{eq:wave-cg-dg-at-once},~\eqref{eq:heat-cg-cg-at-once}
and~\eqref{eq:wave-cg-cg-at-once}. By batching multiple time steps into one, the
implementation allows for $N_{t}$ to become very large. In this case, we plan to
extend our implementation to be able to store the temporal matrices as sparse
matrices. Here, we keep $N_{t}$ at moderate sizes, where this does not pay off
yet.

The implementation of the STMG method is based on the geometric multigrid method
implemented in {\ttfamily
  deal.II}~\cite{munchEfficientDistributedMatrixfree2023}. The only remaining
task is to define the transfer and smoothing operators.

\subsection{Grid transfer}
\label{sec:org8e1deec}
\begin{algorithm}[H]
  \KwResult{Space-time multigrid transfer operators}%
  \For{$\ell\gets0$ \KwTo $L$}{%
    \uIf{mgType[$\ell$] == space}{%
      setup space transfers and transfer operators $\symcal R_{\ell}$,
      $\symcal P_{\ell}$\;%
    }\uElse{%
      setup time transfer operators  $\symcal R_{\ell}$,
      $\symcal P_{\ell}$\;%
      \uIf{mgType[$\ell$] == \(k\)-multigrid}{ $r=r-1$\; }
      \uElseIf{mgType[$\ell$] == \(\tau\)-multigrid}{%
        $time steps /= 2$\;%
      }%
    }%
  }%
  \caption{Algorithm to construct space-time multigrid transfers in a multilevel
    system.}\label{alg:mg-sequence}
\end{algorithm}
The definition of the prolongation operator depends on the coarsening strategy
employed at each grid level, as specified by Algorithm~\ref{alg:mg-sequence}.
For temporal and spatial coarsening, the standard \(L^2\)-projection is applied
by the prolongation operators for $p$ or $k$ and $h$ or $\tau$-multigrid
methods. For a detailed description and performance engineering of the spatial
transfers we refer to~\cite{munchEfficientDistributedMatrixfree2023}. For the
temporal transfers, in the context of discontinuous Galerkin time
discretization, implementing the projections is straightforward. Continuous time
discretizations pose a challenge due to the continuity constraints imposed due
to the Galerkin-Petrov method. In particular, the first temporal degree of
freedom $\vct u_n^0$ for each local subproblem on interval \(I_n\) is determined
by the continuity constraint $\vct u_n^0=\vct u_{n-1}^{N_t}$
(cf.~Section~\ref{sec:orgc1f0a88}) and is therefore considered a constrained
degree of freedom. As a result, the residuals with respect to this degree of
freedom are zero, allowing us to omit them in the transfers of the residuals in
the multigrid method. Of course $\vct u_n^0$ is accounted for by the problem on
$I_{n-1}$, where $\vct u_{n-1}^{N_t}$ is an unknown in the local linear system.
To efficiently implement the temporal transfers, we exploit the tensor product
structure inherent in the discretization. Then, the spatial prolongation
$\symcal{P}_l^h$ applied to the block vector $\symbf{u}_l$ is represented by the
tensor product $\symcal{P}_l^h \otimes \symbf{u}_l$ and the temporal
prolongation is represented by $\symcal{P}_l^\tau \symbf{u}_l$.

Similarly, the choice of the restriction operator is influenced by the selected
coarsening strategy. Typically, the restriction operator is defined as the
adjoint of the prolongation operator. Our implementation provides an option to
compute the restrictions using \(L^{2}\)-projections. However, this approach
negates the desirable property of maintaining adjoint symmetry between the
operators. In practice, both methods have demonstrated comparable performance.
Therefore, we maintain the definition of the restriction operator as the adjoint
of the prolongation operator throughout this paper. Analogous to the
prolongation, the spatial restriction $\symcal R_l^{h}$ is given by
$\symcal R_l^{h} \otimes \symbf{u}_l= {(\symcal P_l^{h})}^{\top} \otimes
\symbf{u}_l$ and the temporal restriction $\symcal R_l^{\tau}$ is given by
$\symcal R_l^{\tau} \symbf{u}_l={(\symcal P_l^{\tau})}^{\top} \symbf{u}_l$.

\subsection{Space-time additive Schwarz smoother}
\label{sec:org8a7ec70}
The smoothing that preceded the restrictions is usually a simple iteration that
approximates the solution of the linear system, i.\,e.\
\[ {\mathtt{smoother}}_{l}(\symbfcal S_l,\,\symbf f_l) \approx \symbfcal S_l^{-1}\symbf
  f_l , \] and aims to quickly reduce all high frequency components of the
residual \(\symbf f_l-\symbfcal S_l\symbf u_l\). If the residual is reduced at a
rate that is constant across all levels \(l\), then the geometric multigrid
method achieves the optimal complexity \(O(N)\). In our implementation we use a
space-time cell wise ASM smoother. The blocks of the ASM consist of all degrees
of freedom of a space-time element, which are solved with an inner direct
solver. The block size of the ASM is then $k(p+1)^d$ for the \(CG(k)-CGP(k)\)
method and $(k+1)(p+1)^d$ for the \(CG(k)-DG(k)\) method. The idea of the
smoother is to solve small subproblems directly and to replace the inverse of
\(\symbfcal S_l^{-1}\) by
\begin{equation}\label{vanka0} {\mathtt{smoother}}_l(\symbfcal S_l,\,\symbf f_l)
  =\left(\sum_{T\in\symcal T_{\tau,\,h}^l}\symbf R_T^T [\symbf R_T \symbfcal S_l\symbf
  R_T^T]^{-1}\symbf R_T\right) \symbf
  f_l
\end{equation}
where \(\symbf R_T\) is the restriction to those nodes that belong to a
space-time mesh element \(T\in\symcal T_{\tau,\,h}^l\) and \(\symbf R_T^T\) is
its transpose. This smoother is computationally expensive and the application on
each cell is of complexity
\begin{equation}
  \label{eq:asm-complex}
O(k^2(p+1)^{2d}) \quad\text{and}\quad
O((k+1)^2(p+1)^{2d})\,,
\end{equation}
for $CG(k)$-$CGP(k)$ and $CG(k)$-$DG(k)$ discretizations. With a relaxation
parameter \(\omega_l\in (0,1)\), an iteration of the smoother (cf.~Algorithm
\ref{alg:stmg}, line~\ref{alg:presmooth} and~\ref{alg:postsmooth}) is then given
by
\begin{equation}\label{vanka}
  \mathtt{Smoother}(\symbfcal S_l,\symbf f_l,\symbf u_l)=\symbf u_l + \omega_l \,{\mathtt{smoother}}_l(\symbfcal S_l,\,\symbf f_l-\symbf S_l\symbf u_l).
\end{equation}
We note that our current implementation of the space-time ASM
smoother does only partially take advantage of the block structure arising from
the tensor-product space-time finite elements. We only exploit the block
structure to prevent the assembly of the full system matrix. Future work will
focus on leveraging these block patterns to improve smoother efficiency and
solver scalability.

\section{Numerical experiments}
\label{sec:org6436fcc}
We evaluate the STMG method when applied to parabolic and hyperbolic PDEs. The
evaluation include numerical convergence tests and two examples of inhomogeneous
coefficients. The investigation focuses on the performance and parallel scaling
of the algebraic solver, as well as its robustness with respect to
discretization order, mesh refinement, and mesh structure. The tests are run on
an HPC cluster (HSUper at HSU) with 571 nodes, each with 2 Intel Xeon Platinum
8360Y CPUs and \SI{256}{\giga\byte} RAM.\@ The processors have 36 cores each.

For the solution of the arising algebraic systems we use the STMG method
developed in this paper as a preconditioner with a GMRES method. The STMG method
(cf. Section \ref{sec:org1919a00}) is employed with a single V-cycle for every
GMRES iteration. In the STMG method we use the space-time cell-based ASM
smoother presented in Section~\ref{sec:org8a7ec70}. The implementation is based
on the geometric multigrid method implemented in {\ttfamily
  deal.II}~\cite{munchEfficientDistributedMatrixfree2023}. Evaluations of the
linear operator $\symbf S$ are performed using the matrix-free infrastructure
described in~\cite{kronbichlerGenericInterfaceParallel2012}. In all numerical
experiments, the stopping criterion for the GMRES iterations is an absolute
residual smaller than $\num{1.e-12}$ (Section~\ref{sec:org177dbb3}) or
$\num{1.e-10}$ (Section~\ref{sec:org1d027f0}) or a reduction of the initial residual by a
factor of $\num{1.e-12}$. The numerical methods are only tested for the $p=k$
case, i.\,e.\ we consider the \(CG(k)\)-\(CGP(k)\) and the \(CG(k)\)-\(DG(k)\)
method. The relaxation parameter $\omega_l$ in~\eqref{vanka} is determined based
on the minimal and maximal eigenvalues $\lambda_{l,\,\min}$,
$\lambda_{l,\,\max}$ on each level by
\begin{equation}
  \label{eq:relaxation}
\omega_l=\frac{2}{\lambda_{l,\,\max} + \lambda_{l,\,\min}}\,.
\end{equation}
The eigenvalues are estimated according to~\cite[Section
2.2]{munchCacheoptimizedLowoverheadImplementations2024}. The relaxation
parameter is optimal for iterative methods applied to symmetric positive
definite matrices, such as the Jacobi method. We adopt it here as a sensible
choice. In the numerical experiments it can lead to improved convergence rates
and performs at least as well as \(\omega_l = 1\).

\subsection{Numerical convergence tests}
\label{sec:org177dbb3}
We verify the space-time finite element methods for the heat and wave equation
for the $CG(k)$-$CGP(k)$ and the $CG(k)$-$DG(k)$ method. In order to be
able to validate the numerical methods by calculating the experimental orders of
convergence (EOC), we prescribe the same solution
\begin{equation}
  \label{eq:exact-prescribed}
  u(\symbfit{x},\,t)=\sin(2\pi f t)\sin(2\pi f x)\sin(2\pi f y)\sin(2\pi f z)\,,
\end{equation}
with $f=2$ for the heat and the wave equation. We consider this 3D test case in
the space-time domain $\Omega\times I= {[0,\,1]}^3\times [0,\,1]$. The initial
time mesh $\symcal{T}_{\tau}$ consists of the twice uniformly refined interval
$I$, while the initial spatial mesh $\symcal{T}_h$ consists of the once
uniformly refined hypercube $\Omega$. The thermal diffusivity and sound speed
are set to 1. All convergence tests were performed on 96 nodes of the HPC
cluster. In all tests we perform a STMG method with at least one level in space
and time. The convergence rate is not affected by the number of time steps that
are batched into a single system. We therefore don't report the convergence
rates for different space-time multigrid sequences. We report the influence of
different space-time multigrid sequences on the performance of the iteration
count. The number of smoothing steps is set to $n_{\text{smooth}}=1$.

\subsubsection{Convergence tests of the heat equation}
\label{sec:convergence-heat}
\begin{figure}[htb]
  \centering
  \includegraphics[width=\linewidth]{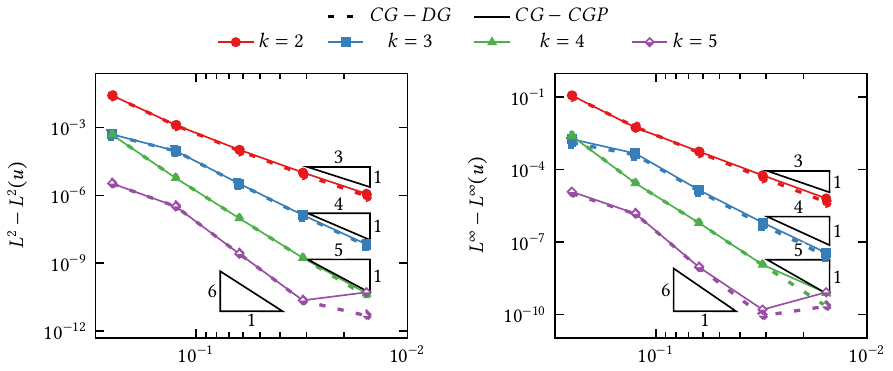}
  \includegraphics[width=\linewidth]{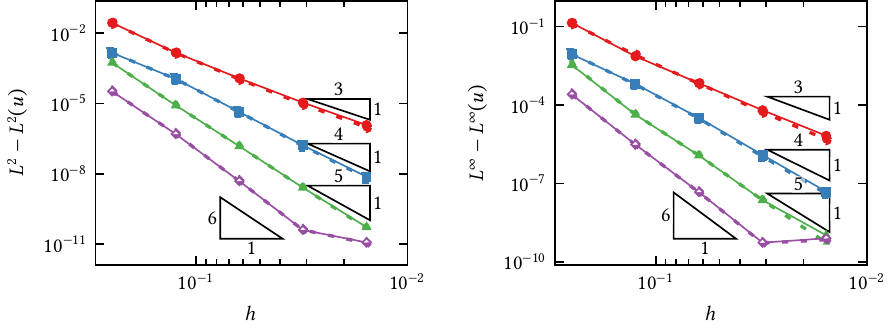}    
  \caption{\label{fig:conv-heat}Calculated errors for $u$ for different orders
    of convergence on Cartesian meshes (top row) and perturbed meshes (bottom
    row) for \(CG(k)-DG(k)\) discretizations of the heat equation. The expected
    orders of convergence $k+1$, represented by the triangles, match with the
    experimental orders.}
\end{figure}

We perform convergence tests for the systems~\eqref{eq:heat-cg-dg-at-once}
and~\eqref{eq:heat-cg-cg-at-once}. We evaluate the performance of STMG on a
Cartesian grid. To ensure effectiveness of STMG on complex spatial domains, we
also consider perturbed meshes. For this, we apply a random displacement to all
vertices of the grid. The direction of movement of each vertex is randomly
generated, while the length of the shift vector is set to 0.15 times the minimal
length of the edges adjacent to this vertex. The convergence test has
homogeneous Dirichlet boundary conditions. We study the errors
\(e_{u}=u(\symbfit x,\,t)-u_{\tau,\,h}(\symbfit x,\,t)\) in the norms given by
\begin{equation}
  \label{eq:norms}
  \norm{e_{u}}_{L^{\infty}(L^{\infty})}=\max_{t\in I}\paran{\sup_{\Omega}\norm{u}_{\infty}}\,,\quad
  \norm{e_{u}}_{L^2(L^2)}=\paran{\int_{I}\int_{\Omega} \lvert
    e_{u}\rvert^{2}\drv \symbfit x\drv t}^{\frac{1}{2}}\,.
\end{equation}
We abbreviate the error quantity \(\norm{e_{u}}_{L^{\infty}(L^2)}\) by
\(L^{\infty}-L^{\infty}(u)\) and analogously for the other norms. We verify the
accuracy of the \(CG(k)-DG(k)\) and \(CG(k)-CGP(k)\) methods for
$p=k,\:k \in \{2,\,3,\,4,\,5\}$. In Figure~\ref{fig:conv-heat} we observe that
the optimal orders of convergence are achieved for Cartesian and perturbed
meshes for \(CG(k)-DG(k)\) discretizations.

We verify the accuracy of the \(CG(k)-DG(k)\) and \(CG(k)-CGP(k)\) methods for
\( p=k, \: k \in \{2, 3, 4, 5\} \). In Figure~\ref{fig:conv-heat} we observe that the optimal orders of
convergence are achieved for Cartesian and perturbed meshes for \(CG(k)-DG(k)\)
discretizations. Therefore, the proposed method is robust and accurate across
different mesh configurations.

Table~\ref{tab:iter-heat} summarizes the average number of GMRES iterations
required per time step to solve the resulting linear systems. The data shows the
solver's grid independence and a moderate increase of the number of iterations
with the order of the discretization. A potential solution to this issue is to
consider $p$-multigrid, as demonstrated in~\cite{fehnHybridMultigridMethods2020}
or by vertex-patch based
smoothers~\cite{pavarinoAdditiveSchwarzMethods1993,miraiMultilevelAlgebraicError2020,schberlAdditiveSchwarzPreconditioning2008}.
On perturbed meshes, there is a slight increase in the number of GMRES
iterations compared to the Cartesian meshes. However, this difference diminishes
with further mesh refinement. Collecting multiple time steps in one system does
not significantly increase the iteration counts.

\begin{figure}[htb]
  \centering
  \includegraphics[width=\linewidth]{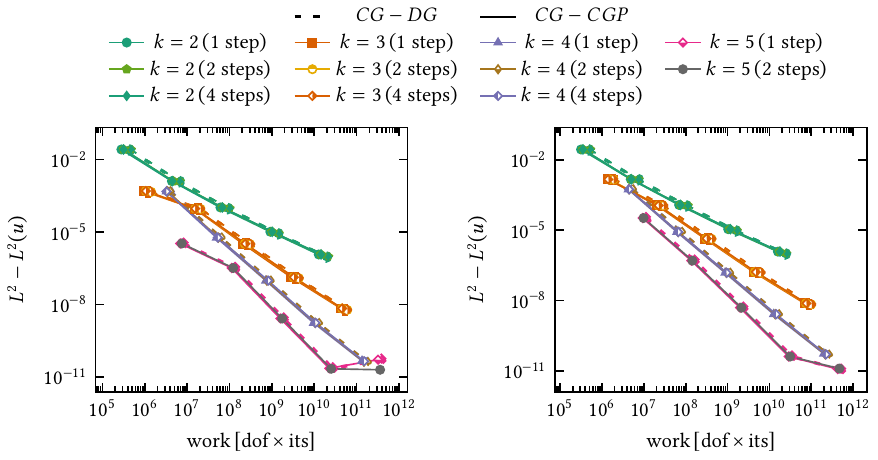}
  \caption{\label{fig:accuracy-work-heat}Calculated $L^2-L^2(u)$-errors of the heat equation for
    different polynomial orders plotted over the work~\eqref{eq:work} on
    Cartesian meshes (left) and perturbed meshes (right). The advantage of
    higher order discretizations can be observed.}
\end{figure}

To quantify the efficiency and investigate the missing robustness with respect
to the polynomial order, Figure~\ref{fig:accuracy-work-heat} presents the
accuracy in terms of the $L^2-L^2$ error over the work \( w \) defined as:
\begin{equation}
  \label{eq:work}
  w = \sum_{i=1}^{\abs{\mathcal{T}_{\tau}}} N^t \times N^{\symbfit{x}} \times N^{it}_i
\end{equation}
where \( N^{it}_i \) denotes the number of GMRES iterations performed for
solving the local subproblem on interval \( I_i \). Thereby, the work accounts
for the computational workload of the global space-time problems.
Figure~\ref{fig:accuracy-work-heat} shows that higher order discretizations are
more efficient. Despite a moderate increase in the number of iterations, they
provide significantly higher accuracy per unit of work. This gain in efficiency,
which balances computational effort with accuracy improvements, proves the
advantages of higher order discretizations.

\begin{table}[ht]
  \caption{\label{tab:iter-heat} Table with the number of GMRES
    iterations until convergence for different polynomial degrees $k$ and number
    of refinements $r$. We show the numbers for the Cartesian mesh (left) and
    the perturbed mesh (right) for the \(CG(k)-DG(k)\) discretization of the
    heat equation.}
  \centering
  \begin{subcaptionblock}{\textwidth}
    \centering
    \caption{\(CG(k)-DG(k)\) discretizations on a Cartesian mesh
      (left) and a perturbed mesh (right).}\label{tab:iter-heat-dg}
    \begin{minipage}{0.45\textwidth}
      \begin{center}\small
        \begin{tabular}{crrrrr}
          \toprule
          \(k\) $\backslash$ \(r\)&    2&       3&       4&       5&       6\\
          \midrule
          2&  9.0&  8.75&  7.75&  6.813& 5.938\\
          3& 11.0&  9.75&  8.88&  7.688& 6.188\\
          4& 14.0& 12.75& 10.88&  9.500& 7.688\\
          5& 16.5& 14.25& 12.75& 10.688& 8.568\\
          \bottomrule
        \end{tabular}
      \end{center}
    \end{minipage}
    \hspace{.5cm}
    \begin{minipage}{0.45\textwidth}
      \begin{center}\small
        \begin{tabular}{crrrrr}
          \toprule
          \(k\) $\backslash$ \(r\)&    2&       3&       4&       5&       6\\
          \midrule
          2&  9.0&  9.75&  9.00&  8.875&  8.656\\
          3& 12.0& 11.75& 10.88& 10.188& 10.563\\
          4& 14.5& 14.00& 12.88& 11.813& 11.781\\
          5& 15.0& 14.50& 13.63& 12.563& 11.781\\
          \bottomrule
        \end{tabular}
      \end{center}
    \end{minipage}
  \end{subcaptionblock}%

  \begin{subcaptionblock}{\textwidth}
    \centering
    \caption{\(CG(k)-CGP(k)\) discretizations on a Cartesian mesh
      (left) and a perturbed mesh (right).}\label{tab:iter-heat-cgp}
    \begin{minipage}{0.45\textwidth}
      \begin{center}\small
        \begin{tabular}{crrrrr}
          \toprule
          \(k\) $\backslash$ \(r\)&    2&       3&       4&       5&       6\\
          \midrule
          2&  8&  8.75&  8.125& 7.875& 6.938\\
          3&  8&  8.75&  7.875& 6.938& 6.563\\
          4& 11& 11.25& 10.625& 9.063& 7.938\\
          5& 11& 12.00& 10.875& 9.875& 8.781\\
          \bottomrule
        \end{tabular}
      \end{center}
    \end{minipage}
    \hspace{.5cm}
    \begin{minipage}{0.45\textwidth}
      \begin{center}\small
        \begin{tabular}{crrrrr}
          \toprule
          \(k\) $\backslash$ \(r\)&    2&       3&       4&       5&       6\\
          \midrule
           2&  9.0&  9.75&  9.25&  8.875&  8.688\\
           3& 12.0& 12.00& 10.88& 10.188& 10.594\\
           4& 14.5& 14.00& 12.88& 11.875& 11.781\\
           5& 15.0& 14.50& 13.63& 12.563& 11.813\\
          \bottomrule
        \end{tabular}
      \end{center}
    \end{minipage}
  \end{subcaptionblock}%
  
  \begin{subcaptionblock}{\textwidth}
    \centering
    \caption{\(CG(k)-DG(k)\) with 2 time steps at
      once on a Cartesian mesh
      (left) and a perturbed mesh (right).}\label{tab:iter-heat-dg-2s}
    \begin{minipage}{0.45\textwidth}
      \begin{center}\small
        \begin{tabular}{crrrrr}
          \toprule
          \(k\) $\backslash$ \(r\)&    2&       3&       4&       5&       6\\
          \midrule
          2&  9&  9&  9.00&    8.344& 7.813\\
          3&  9&  9&  8.69&    7.875& 6.938\\
          4& 11& 12& 11.19&    9.875& 8.813\\
          \bottomrule
        \end{tabular}
      \end{center}
    \end{minipage}
    \hspace{.5cm}
    \begin{minipage}{0.45\textwidth}
      \begin{center}\small
        \begin{tabular}{crrrrr}
          \toprule
          \(k\) $\backslash$ \(r\)&    2&       3&       4&       5&       6\\
          \midrule
          2& 10& 10.0&  10.00&  9.60&  9.234\\
          3& 12& 12.38& 11.75& 10.88& 11.484\\
          4& 15& 15.0&  13.75& 12.88& 12.750\\
          \bottomrule
        \end{tabular}
      \end{center}
    \end{minipage}
  \end{subcaptionblock}%

  \begin{subcaptionblock}{\textwidth}
    \centering
    \caption{\(CG(k)-CGP(k)\) with 2 time steps at
      once on a Cartesian mesh
      (left) and a perturbed mesh (right).}\label{tab:iter-heat-cgp-2s}
    \begin{minipage}{0.45\textwidth}
      \begin{center}\small
        \begin{tabular}{crrrrr}
          \toprule
          \(k\) $\backslash$ \(r\)&    2&       3&       4&       5&       6\\
          \midrule
          2&  9.00&  9.38&  9.00& 8.59& 7.813\\
          3&  9.00&  9.00&  8.68& 7.88& 6.938\\
          4& 11.75& 12.38& 11.19& 9.88& 8.813\\
          \bottomrule
        \end{tabular}
      \end{center}
    \end{minipage}
    \hspace{.5cm}
    \begin{minipage}{0.45\textwidth}
      \begin{center}\small
        \begin{tabular}{crrrrr}
          \toprule
          \(k\) $\backslash$ \(r\)&    2&       3&       4&       5&       6\\
          \midrule
           2& 10.0& 10.00& 10.00&  9.750&  9.484\\
           3& 12.8& 13.00& 11.75& 10.875& 11.484\\
           4& 15.0& 15.00& 13.75& 12.875& 12.750\\
          \bottomrule
        \end{tabular}
      \end{center}
    \end{minipage}
  \end{subcaptionblock}%
  
  \begin{subcaptionblock}{\textwidth}
    \centering
    \caption{\(CG(k)-DG(k)\) with 4 time steps at
      once on a Cartesian mesh
      (left) and a perturbed mesh (right).}\label{tab:iter-heat-dg-4s}
    \begin{minipage}{0.45\textwidth}
      \begin{center}\small
        \begin{tabular}{crrrrr}
          \toprule
          \(k\) $\backslash$ \(r\)&    2&       3&       4&       5&       6\\
          \midrule
          2&  9&  9&  9.0&    8.5& 7.75\\
          3&  9&  9&  9.0&    8.0& 7.00\\
          4& 11& 12& 11.5&   10.0& 9.50\\
          \bottomrule
        \end{tabular}
      \end{center}
    \end{minipage}
    \hspace{.5cm}
    \begin{minipage}{0.45\textwidth}
      \begin{center}\small
        \begin{tabular}{crrrrr}
          \toprule
          \(k\) $\backslash$ \(r\)&    2&       3&       4&       5&       6\\
          \midrule
          2&  10& 10&  10&  9.75&  9.625\\
          3&  12& 13&  12& 11.50& 11.750\\
          4&  15& 15&  14& 13.25& 13.125\\
          \bottomrule
        \end{tabular}
      \end{center}
    \end{minipage}
  \end{subcaptionblock}%

  \begin{subcaptionblock}{\textwidth}
    \centering
    \caption{\(CG(k)-CGP(k)\) with 4 time steps at
      once on a Cartesian mesh
      (left) and a perturbed mesh (right).}\label{tab:iter-heat-cgp-4s}
    \begin{minipage}{0.45\textwidth}
      \begin{center}\small
        \begin{tabular}{crrrrr}
          \toprule
          \(k\) $\backslash$ \(r\)&    2&       3&       4&       5&       6\\
          \midrule
          2&  9.00 &  9.00&  9.00& 8.75& 7.88\\
          3&  9.00 &  9.00&  9.00& 8.00& 7.00\\
          4& 11.00 & 12.00& 12.00&10.00& 8.88\\
          \bottomrule
        \end{tabular}
      \end{center}
    \end{minipage}
    \hspace{.5cm}
    \begin{minipage}{0.45\textwidth}
      \begin{center}\small
        \begin{tabular}{crrrrr}
          \toprule
          \(k\) $\backslash$ \(r\)&    2&       3&       4&       5&       6\\
          \midrule
           2& 10&  10& 10&  10.0&  9.75\\
           3& 12&  13& 12&  11.5& 11.75\\
           4& 15&  15& 14&  13.5& 13.25\\
          \bottomrule
        \end{tabular}
      \end{center}
    \end{minipage}
  \end{subcaptionblock}%
\end{table}
\FloatBarrier%
\subsubsection{Convergence tests of the wave equation}
\label{sec:org9bc35f4}
\begin{figure}[htbp]
  \centering
  \includegraphics[width=\linewidth]{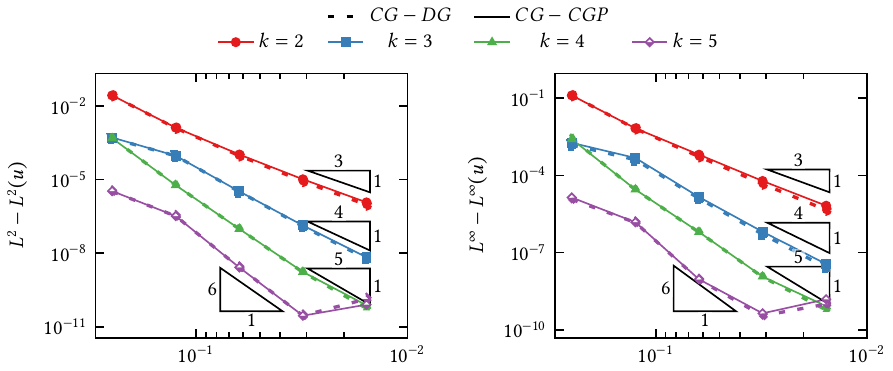}
  \includegraphics[width=\linewidth]{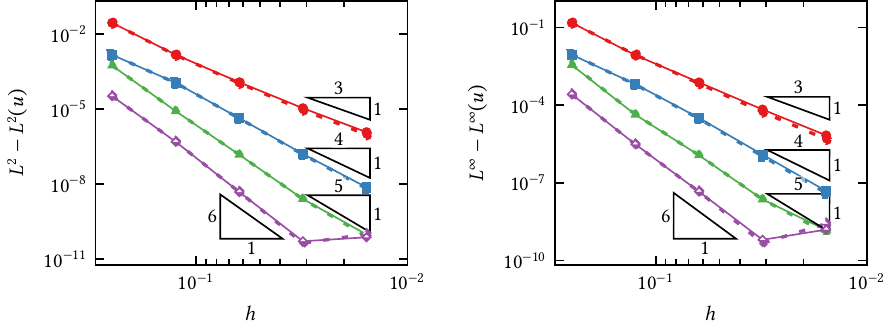}    
  \caption{\label{fig:conv-wave}Calculated errors for the displacement $u$ for
    different orders on Cartesian (top row) and perturbed meshes (bottom row)
    for \(CG(k)-DG(k)\) (dashed) and \(CG(k)-CGP(k)\) (solid) discretizations of
    the wave equation. The expected orders of convergence $k+1$, match with the
    experimental orders.}
\end{figure}
\begin{figure}[htbp]
  \includegraphics[width=\linewidth]{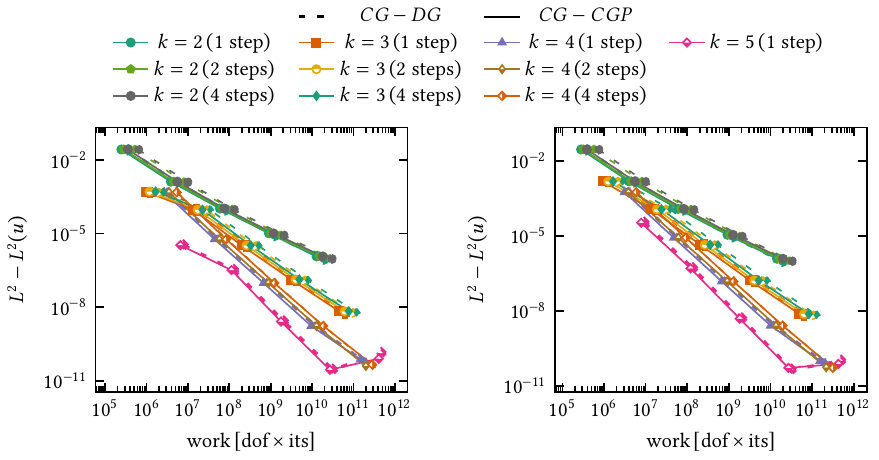}    
  \caption{\label{fig:accuracy-work-wave}Calculated $L^2-L^2(u)$-errors of the
    wave equation for different orders plotted over the work on
    Cartesian meshes (left) and perturbed meshes (right). Higher-order
    discretizations are beneficial.}
\end{figure}
We perform the convergence tests for the systems~\eqref{eq:wave-cg-dg-at-once}
and~\eqref{eq:wave-cg-cg-at-once}. Similar to the heat equation in
Section~\ref{sec:convergence-heat}, we test the wave equation on a Cartesian and
a perturbed mesh. The perturbed mesh is generated analogously to the one in
Section~\ref{sec:convergence-heat}.

We verify the accuracy of the \(CG(k)-DG(k)\) and \(CG(k)-CGP(k)\) methods for
\( p=k, \: k \in \{2, 3, 4, 5\} \). In Figure~\ref{fig:conv-wave} we observe that the optimal orders of
convergence are achieved for Cartesian and perturbed meshes for \(CG(k)-DG(k)\)
and \(CG(k)-CGP(k)\) discretizations. The proposed method is robust and accurate
across different mesh configurations for the heat and wave equation. The
continuous and discontinuous time discretizations lead to similar errors.

Table~\ref{tab:iter-wave} summarizes the average number of GMRES iterations
required per time step to solve the resulting linear systems. Similar to the
heat equation, the data shows the solver's grid independence, but we don't
achieve full independence of the order of the discretization. On perturbed
meshes, there is a slight increase in the number of GMRES iterations compared to
the Cartesian meshes. However, this difference vanishes with the mesh
refinement. With regard to the iteration counts, the \(CG(k)-CGP(k)\) and
\(CG(k)-DG(k)\) discretizations in Table~\ref{tab:iter-wave} exhibit comparable
performance. The solver's performance is similar for the wave and heat equation,
except for the fact, that the collection of multiple time steps in one system
leads to an increase of the number of GMRES iterations.

To quantify the efficiency and investigate the missing robustness with respect
to the polynomial order, Figure~\ref{fig:accuracy-work-wave} presents the
accuracy in terms of the $L^2-L^2$ error over the work \( w \) as defined
in~\eqref{eq:work}. Analogous to the investigations on the heat equation,
Figure~\ref{fig:accuracy-work-wave} shows that higher order discretizations are
more efficient. Despite a moderate increase in the number of iterations, they
provide significantly higher accuracy per unit of work. This gain in efficiency
can be observed for the heat and wave equation in the context of the proposed
STMG method. In Figure~\ref{fig:accuracy-work-wave} we observe
that, in contrast to Figure~\ref{fig:accuracy-work-heat}, the amount of work
increases with the number of time steps collected in one system. This is of
course expected due to the increase of GMRES iterations. Further,
Figure~\ref{fig:accuracy-work-wave} illustrates that the continuous
discretization method exhibits an advantage over the discontinuous time
discretizations. This is due to the smaller number of temporal degrees of
freedom resulting from the continuity constraint, as discussed in
Section~\ref{sec:disc-cgp}. Due to the large scales involved, the differences
are understated in the plot. To illustrate, consider the case of k = 5 with the
largest refinement, resulting in 5 and 6 temporal degrees of freedom per local
subproblem on $I_n$ for the continuous and discontinuous time discretization.
Therefore, the number of global degrees of freedom is $\num{36522640000}$ and
$\num{43827168000}$ for \(CG(k)-CGP(k)\) and \(CG(k)-DG(k)\). Whether the
advantage of having $\frac{1}{k+1}$ less degrees of freedom for the same level
of accuracy materializes in a decreased time to solution depends on the number
of iterations (cf.~Table~\ref{tab:iter-wave}), the underlying hardware and
scaling of the software and on the problem. In the context of dynamic
poroelasticity simulations for instance, the authors
of~\cite{anselmannBenchmarkComputationsDynamic2023} observed that continuous
time discretizations are efficient, yet not clearly superior to discontinuous
methods. In the convergence tests within this section, \(CG(k)-DG(k)\)
discretizations exhibited runtimes nearly double those of their continuous
counterparts. This unexpected difference will be investigated in the following
section by extensive scaling tests. The node configurations and simulation
parameters in this section were selected to ensure that the simulations are
completed within a specified timeframe and not optimized for optimal
performance.

\begin{table}[htb]
  \centering
  \caption{Tables with the number of GMRES iterations until convergence for
    different polynomial degrees $k$ and number of refinements $r$. We show the
    numbers for the Cartesian mesh (left) and the perturbed mesh (right) for the
    \(CG(k)-DG(k)\) discretization of the wave equation.}\label{tab:iter-wave}
  \begin{subcaptionblock}{\textwidth}
    \centering
    \caption{\(CG(k)-DG(k)\) discretizations on a Cartesian mesh
      (left) and a perturbed mesh (right).}\label{tab:iter-wave-dg}
    \begin{minipage}{0.45\textwidth}
      \begin{center}\small
        \begin{tabular}{crrrrr}
          \toprule
          \(k\) $\backslash$ \(r\)&    2&       3&       4&       5&       6\\
          \midrule
          2&  6.5&  7.00&  7.000&  6.938&  6.906\\
          3&  8.0&  7.75&  7.875&  7.813&  7.563\\
          4&  9.0&  9.75&  9.875&  9.375&  8.813\\
          5& 10.5& 11.50& 11.875& 11.563& 10.844\\
          \bottomrule
        \end{tabular}
      \end{center}
    \end{minipage}
    \hspace{.5cm}
    \begin{minipage}{0.45\textwidth}
      \begin{center}\small
        \begin{tabular}{crrrrr}
          \toprule
          \(k\) $\backslash$ \(r\)&    2&       3&       4&       5&       6\\
          \midrule
          2&  7.0&  7.00&  7.000&  6.938&  6.906\\
          3&  8.5&  8.00&  8.000&  7.938&  7.844\\
          4& 10.5& 10.75&  9.875&  9.500&  9.156\\
          5& 12.5& 12.25& 11.875& 11.625& 11.375\\
          \bottomrule
        \end{tabular}
      \end{center}
    \end{minipage}
  \end{subcaptionblock}%

  \begin{subcaptionblock}{\textwidth}
    \centering
    \caption{\(CG(k)-CGP(k)\) discretizations on a Cartesian mesh
      (left) and a perturbed mesh (right).}\label{tab:iter-wave-cgp}
    \begin{minipage}{0.45\textwidth}
      \begin{center}\small
        \begin{tabular}{crrrrr}
          \toprule
          \(k\) $\backslash$ \(r\)&    2&       3&       4&       5&       6\\
          \midrule
          2&  7.0&  7.75&  7.875&  7.750&  6.969\\
          3&  8.0&  7.75&  7.625&  7.625&  6.969\\
          4&  9.0& 10.00&  9.875&  8.938&  8.750\\
          5& 10.5& 12.00& 12.125& 11.750& 11.219\\
          \bottomrule
        \end{tabular}
      \end{center}
    \end{minipage}
    \hspace{.5cm}
    \begin{minipage}{0.45\textwidth}
      \begin{center}\small
        \begin{tabular}{crrrrr}
          \toprule
          \(k\) $\backslash$ \(r\)&    2&       3&       4&       5&       6\\
          \midrule
          2&  8.0&  8.00&  8.000&  7.938&  7.813\\
          3&  8.0&  8.00&  7.875&  7.813&  7.625\\
          4& 10.5& 10.75&  9.875&  9.438&  9.219\\
          5& 12.5& 12.75& 12.500& 11.750& 11.750\\
          \bottomrule
        \end{tabular}
      \end{center}
    \end{minipage}
  \end{subcaptionblock}%
  
  \begin{subcaptionblock}{\textwidth}
    \centering
    \caption{\(CG(k)-DG(k)\) with 2 time steps at
      once on a Cartesian mesh
      (left) and a perturbed mesh (right).}\label{tab:iter-wave-dg-2s}
    \begin{minipage}{0.45\textwidth}
      \begin{center}\small
        \begin{tabular}{crrrrr}
          \toprule
          \(k\) $\backslash$ \(r\)&    2&       3&       4&       5&       6\\
          \midrule
          2&  10.0& 10.88&  10.44& 10.00&  9.73\\
          3&  12.0& 11.88&  11.94& 11.47& 11.73\\
          4&  14.8& 15.88&  14.94& 14.44& 13.77\\
          \bottomrule
        \end{tabular}
      \end{center}
    \end{minipage}
    \hspace{.5cm}
    \begin{minipage}{0.45\textwidth}
      \begin{center}\small
        \begin{tabular}{crrrrr}
          \toprule
          \(k\) $\backslash$ \(r\)&    2&       3&       4&       5&       6\\
          \midrule
          2&  10& 10.00&  10.00&  9.60&  9.234\\
          3&  12& 12.38&  11.75& 10.88& 11.484\\
          4&  15& 15.00&  13.75& 12.88& 12.750\\
          \bottomrule
        \end{tabular}
      \end{center}
    \end{minipage}
  \end{subcaptionblock}%

  \begin{subcaptionblock}{\textwidth}
    \centering
    \caption{\(CG(k)-CGP(k)\) with 2 time steps at
      once on a Cartesian mesh
      (left) and a perturbed mesh (right).}\label{tab:iter-wave-cgp-2s}
    \begin{minipage}{0.45\textwidth}
      \begin{center}\small
        \begin{tabular}{crrrrr}
          \toprule
          \(k\) $\backslash$ \(r\)&    2&       3&       4&       5&       6\\
          \midrule
          2&  9.00&  9.38&  9.00& 8.59& 7.813\\
          3&  9.00&  9.00&  8.68& 7.88& 6.938\\
          4& 11.75& 12.38& 11.19& 9.88& 8.813\\
          \bottomrule
        \end{tabular}
      \end{center}
    \end{minipage}
    \hspace{.5cm}
    \begin{minipage}{0.45\textwidth}
      \begin{center}\small
        \begin{tabular}{crrrrr}
          \toprule
          \(k\) $\backslash$ \(r\)&    2&       3&       4&       5&       6\\
          \midrule
           2& 10.0& 10& 10.00&  9.750&  9.484\\
           3& 12.8& 13& 11.75& 10.875& 11.484\\
           4& 15.0& 15& 13.75& 12.875& 12.750\\
          \bottomrule
        \end{tabular}
      \end{center}
    \end{minipage}
  \end{subcaptionblock}%

  \begin{subcaptionblock}{\textwidth}
    \centering
    \caption{\(CG(k)-DG(k)\) with 4 time steps at
      once on a Cartesian mesh
      (left) and a perturbed mesh (right).}\label{tab:iter-wave-dg-4s}
    \begin{minipage}{0.45\textwidth}
      \begin{center}\small
        \begin{tabular}{crrrrr}
          \toprule
          \(k\) $\backslash$ \(r\)&    2&       3&       4&       5&       6\\
          \midrule
          2&  12& 13& 11.5& 11.75& 10.88\\
          3&  16& 15& 14.5& 14.50& 13.75\\
          4&  20& 20& 19.0& 18.50& 18.88\\
          \bottomrule
        \end{tabular}
      \end{center}
    \end{minipage}
    \hspace{.5cm}
    \begin{minipage}{0.45\textwidth}
      \begin{center}\small
        \begin{tabular}{crrrrr}
          \toprule
          \(k\) $\backslash$ \(r\)&    2&       3&       4&       5&       6\\
          \midrule
          2& 14& 13&  12.5&  11.75&  11.88\\
          3& 18& 17&  15.5&  15.00&  14.88\\
          4& 22& 22&  20.0&  19.50&  19.00 \\
          \bottomrule
        \end{tabular}
      \end{center}
    \end{minipage}
  \end{subcaptionblock}%

  \begin{subcaptionblock}{\textwidth}
    \centering
    \caption{\(CG(k)-CGP(k)\) with 4 time steps at
      once on a Cartesian mesh
      (left) and a perturbed mesh (right).}\label{tab:iter-wave-cgp-4s}
    \begin{minipage}{0.45\textwidth}
      \begin{center}\small
        \begin{tabular}{crrrrr}
          \toprule
          \(k\) $\backslash$ \(r\)&    2&       3&       4&       5&       6\\
          \midrule
          2&  10&  11&  10.0& 10.0 & 9.88\\
          3&  14&  13&  12.5& 12.0& 11.88\\
          4&  18&  19&  18.0& 17.5& 16.75\\
          \bottomrule
        \end{tabular}
      \end{center}
    \end{minipage}
    \hspace{.5cm}
    \begin{minipage}{0.45\textwidth}
      \begin{center}\small
        \begin{tabular}{crrrrr}
          \toprule
          \(k\) $\backslash$ \(r\)&    2&       3&       4&       5&       6\\
          \midrule
           2&  11 & 11& 11.0& 10.75& 10.88\\
           3&  14 & 14& 13.5& 12.75& 12.88\\
           4&  20 & 21& 19.0& 18.50& 18.75\\
          \bottomrule
        \end{tabular}
      \end{center}
    \end{minipage}
  \end{subcaptionblock}%
\end{table}

\FloatBarrier%
\subsection{Practical example motivated by structural health monitoring}
\label{sec:org1d027f0}
We consider a 3D test problem inspired by
\autocite[Example 5.4]{bangerthAdaptiveGalerkinFinite2010}. Let $\Omega\times I =
{[-1,\,1]}^3\times[0,\,2]\subset \R^3\times \R$ and the initial values be
\begin{equation}
  \label{eq:initial-shm}
  u_0^0(\symbfit{x})=\e^{-\abs{\symbfit x s^{-1}}^2}(1-\abs{\symbfit x
    s^{-1}}^2)\Theta(1-\abs{\symbfit xs^{-1}}),\quad v=0\,,
\end{equation}
where $\Theta$ is the Heaviside function and $s=0.01$. We choose the speed of
sound as
\begin{equation}\label{eq:rho}
\rho =
\begin{cases}
1 & \text{for } y < 0.2, \\
9 & \text{for } y \ge 0.2 \text{ and } z < 0.2, \\
16 & \text{for } y \ge 0.2 \text{ and } z \ge 0.2\,.
\end{cases}
\end{equation}
To resolve the coefficients, we utilize a uniform Cartesian mesh consisting of
$5\times 5\times 5$ cells as an initial coarse mesh. The problem is relevant to
critical applications in structural health monitoring, geophysics, and
seismology. The ability to accurately predict the arrival time of signals at
specific locations is essential for the analysis of structural integrity and
seismic events in these fields. The test problem serves as a foundational model
for simulating complex wave propagation. As goal quantities we define point
evaluations of the displacement $u$ at the 3 points
$\symbfit x_1={(0.75,\,0,\,0)}^{\top}$, $\symbfit x_2={(0, 0, 0.75)}^{\top}$,
$\symbfit x_3={(0.75, 0.1, 0.75)}^{\top}$. The computational tests use up to
$256$ nodes and all cores of the compute nodes are used.

\begin{figure}[htbp]
  \includegraphics[width=\linewidth]{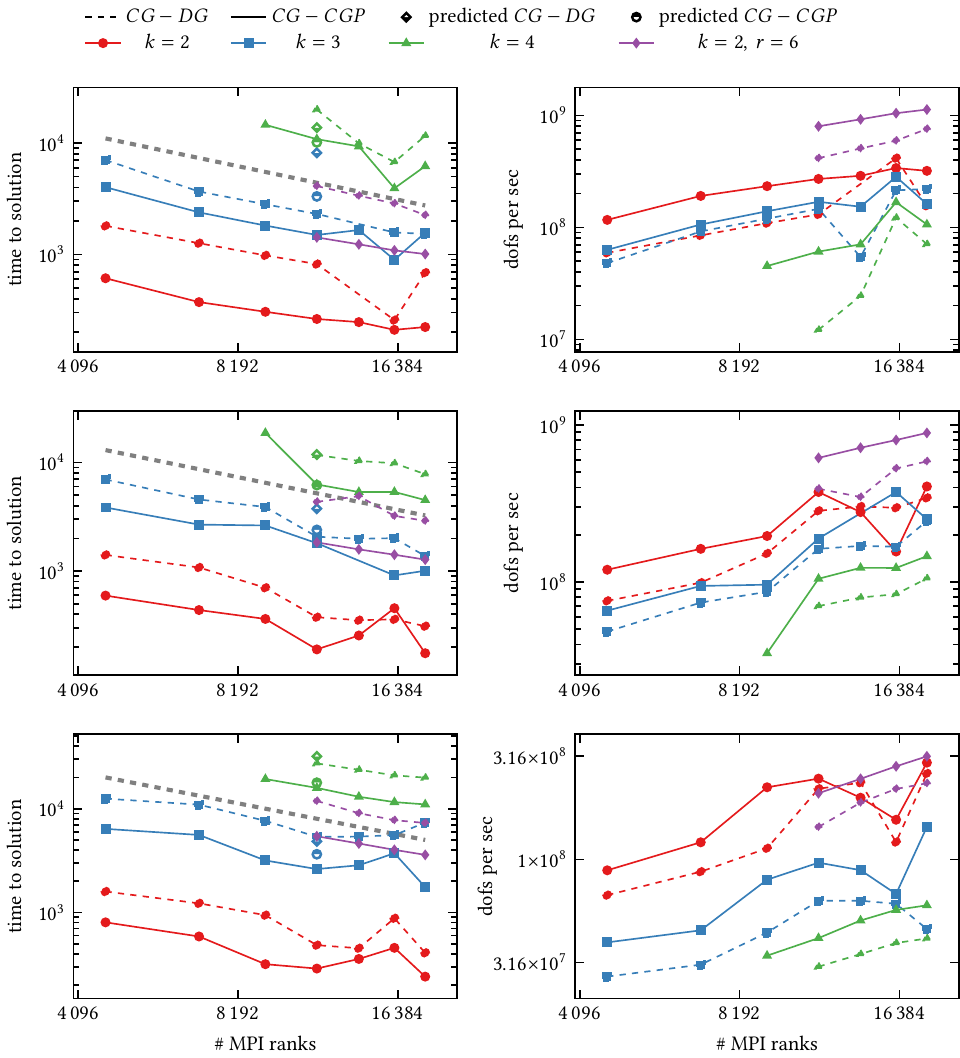}
  \caption{\label{fig:rough-strong-scale}Strong scaling test results for the
    STMG algorithm with varying numbers of smoothing steps. Each row corresponds
    to a different number of smoothing steps: 1, 2, and 4 from top to bottom.
    The left column shows the time to solution as a function of the number of
    MPI processes. The dashed gray lines indicate the optimal scaling. The
    half-squares and half-circles indicate the predicted runtime of the
    \(CG(k)\)-\(DG(k)\) and \(CG(k)\)-\(CGP(k)\) discretization, based on the
    wallclock time of one degree lower and~\eqref{eq:asm-complex}. The right
    column depicts the degrees of freedom (dofs) processed per second over the
    number of MPI processes.}
\end{figure}

\paragraph*{Heterogeneous Coefficients}
Similar to the numerical convergence tests, the initial time mesh
$\symcal{T}_{\tau}$ consists of an $r+1$ times uniformly refined interval $I$,
while the initial spatial mesh $\symcal{T}_h$ consists of the $r$ times
uniformly refined hypercube $\Omega$. For the strong scaling test shown in
Figure~\ref{fig:rough-strong-scale} we set $r=5$, which results in $\num{320}$
time cells and $\num{4096000}$ space cells. We were able to obtain a value $r=6$
and $k=2$. Higher orders for r=6 were not tested due to a lack of computational
resources. The number of spatial, temporal and global degrees of freedom
dofs$_{\Omega}$, dofs$_{\tau}$, dofs$_{\text{global}}$, average number of GMRES
iterations $\overline{n}_{\text{iter}}$ and maximal throughput for this setting
are collected in Table~\ref{tab:iter-wave-hc}.

In our study, we test the STMG algorithm with a fixed number of smoothing steps
$n_{\text{smooth}}\in \{1,\,2,\,4\}$. The variable $V$-cycle proposed
in~\cite{brambleNewConvergenceEstimates1987} was tested but did not lead to an
advantage in terms of iteration count and runtimes. When the variable setting is
used, the number of smoothing steps is doubled at each coarser level. When the
coarse mesh is too small to have cells present on all processors this lead to
significantly decreased parallel efficiency in our experiments.
\begin{table}[htb]
  \centering
  \caption{Tables with the number of GMRES iterations until convergence for
    different polynomial degrees $k$ and number of refinements $r$. We show the
    numbers for the \(CG(k)-DG(k)\) and \(CG(k)-CGP(k)\) discretization of the
    wave equation. The last entry seperated by a rule is a value for $r=6$.}\label{tab:iter-wave-hc}
    \begin{subcaptionblock}{\textwidth}
      \centering
      \caption{Iterations for \(CG(k)-DG(k)\) discretizations with different numbers of smoothing steps ($n_{\text{smooth}}$).}\label{tab:iter-wave-dg-hc}
      \footnotesize
      \begin{tabular}{rrrrrrrrrr}
        \toprule
        $k$ & dofs$_{\Omega}$ &  dofs$_{\tau}$ &  dofs$_{\text{global}}$ & \multicolumn{2}{c}{$n_{\text{smooth}}=1$}                        & \multicolumn{2}{c}{$n_{\text{smooth}}=2$}                   & \multicolumn{2}{c}{$n_{\text{smooth}}=4$}\\
          &                &                &                       & $\overline{n}_{\text{iter}}$ & $\max \text{dofs}/\text{sec}$     & $\overline{n}_{\text{iter}}$& $\max \text{dofs}/\text{sec}$ & $\overline{n}_{\text{iter}}$& $\max \text{dofs}/\text{sec}$\\
        \midrule
        2 & \num{111284641} & $3\times 320$ & \num{106833255360} & 13.73 & \num{414886429} & 8.78  & \num{343515291}& 5.82  & \num{261014550}\\
        3 & \num{263374721} & $4\times 320$ & \num{337119642880} & 13.97 & \num{220339636} & 8.81  & \num{244999740}& 7.32  & \num{62989338}\\
        4 & \num{513922401} & $5\times 320$ & \num{822275841600} & 19.79 & \num{121872808} & 12.32 & \num{105419980}& 10.34 & \num{41326994}\\
        \midrule
        2 & \num{887503681} & $3\times 640$ & \num{1704007067520} & 12.99 & \num{756328037} & 8.87  & \num{586779293}& 5.76  & \num{234034757}\\
        \bottomrule
        \end{tabular}
    \end{subcaptionblock}%

  \begin{subcaptionblock}{\textwidth}
    \centering
    \caption{Iterations for \(CG(k)-CGP(k)\) discretizations with different numbers of smoothing steps.}\label{tab:iter-wave-cgp-hc}
    \footnotesize
    \begin{tabular}{rrrrrrrrrr}
      \toprule
      $k$ & dofs$_{\Omega}$ &  dofs$_{\tau}$ &  dofs$_{\text{global}}$ & \multicolumn{2}{c}{$n_{\text{smooth}}=1$}                        & \multicolumn{2}{c}{$n_{\text{smooth}}=2$}                   & \multicolumn{2}{c}{$n_{\text{smooth}}=4$}\\
        &                &                &                       & $\overline{n}_{\text{iter}}$ & $\max \text{dofs}/\text{sec}$     & $\overline{n}_{\text{iter}}$& $\max \text{dofs}/\text{sec}$ & $\overline{n}_{\text{iter}}$& $\max \text{dofs}/\text{sec}$\\
      \midrule
      2 & \num{111284641} & $2\times 320$ &\num{71222170240} & 11.51 & \num{338508413} & 7.83  &\num{406055703}& 5.46 &\num{294671784} \\
      3 & \num{263374721} & $3\times 320$ &\num{252839732160}& 11.06 & \num{282565637} & 6.88  &\num{251082157} & 5.77 &\num{143822373} \\
      4 & \num{513922401} & $4\times 320$ &\num{657820673280}& 16.25 & \num{167640334} & 10.13 &\num{146020127} & 8.16 &\num{59856294} \\
      \midrule
      2 & \num{887503681} & $2\times 640$ & \num{1136004711680} & 11.19 & \num{1126988801} & 7.82  & \num{890984088}& 5.62  & \num{316259664}\\
      \bottomrule
    \end{tabular}
  \end{subcaptionblock}%
\end{table}

\begin{figure}[htbp]
  \centering
  \includegraphics[width=\linewidth]{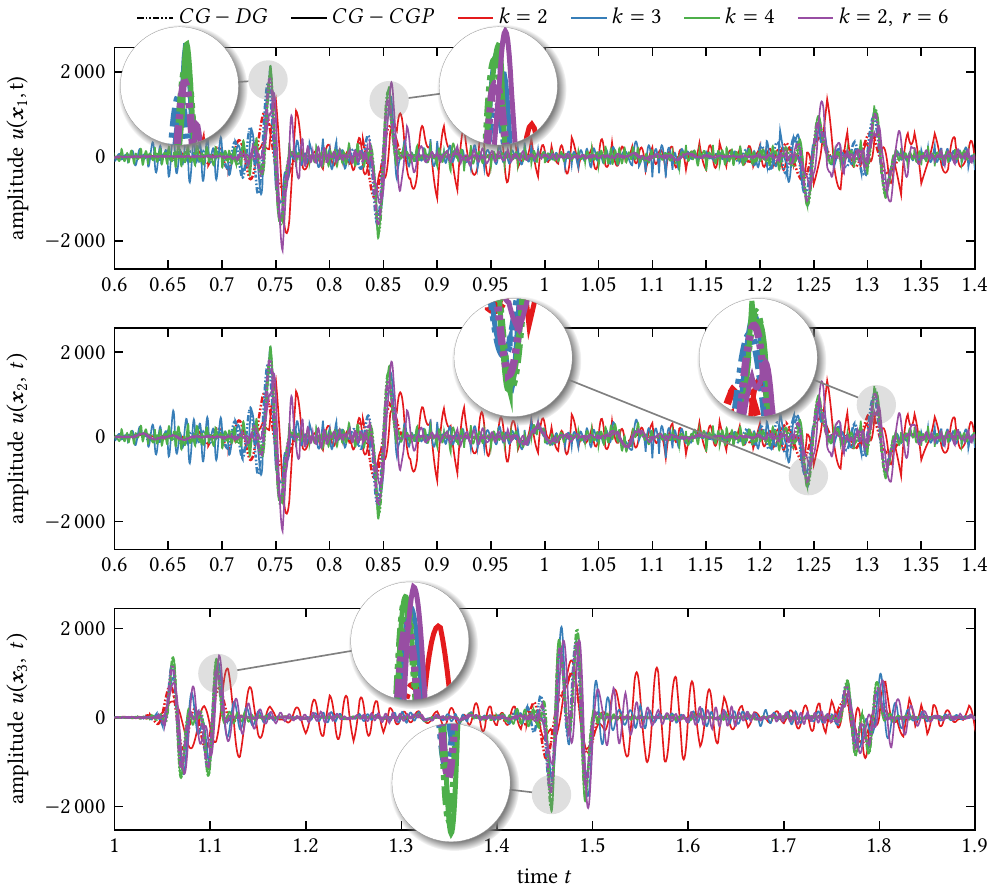}
  \caption{\label{fig:rough-goal-quants}Plots of the point evaluations of the
    displacement at $\symbfit x_1,\:\symbfit x_2,\:\symbfit
    x_3$ plotted over subintervals of the time interval
    $[0,\,2]$. They are adapted such that the first signal arriving at the point
    is plotted.}
\end{figure}

We evaluate the performance of the STMG method. The key ingredient for a
multigrid method is the smoother. The results demonstrate satisfactory
performance of the ASM smoother. However, it is important to note that the
computational cost is significant due to the use of an inner direct solver,
although it is only space-time cell-wise. The cost increases significantly for
higher polynomial orders, which presents a challenge for scalability, as the
computational costs for the smoother begins to dominate. We observe this in
Figure~\ref{fig:rough-strong-scale}, where we plot predictions of the wall times
for discretizations of degree $k$ based on the wall times of the discretizations
of degree $k-1$ and the computational complexity of the
smoother~\ref{eq:asm-complex}. The predictions get increasingly better for more
smoothing steps and higher order. Therefore, the runtime and scaling are
primarily determined by the smoother. Later, we also substantiate this further
for a similar test with larger inhomogeneity. However, the dominance of the
smoother attests to the high performance and efficiency of the matrix-free
framework. We observe in Figure~\ref{fig:rough-strong-scale} that the scaling is
optimal. For both discretizations, the number of GMRES iterations required for
convergence generally increases with higher polynomial degrees \(k\). Increasing
the number of smoothing steps \(n_{\text{smooth}}\) typically reduces the number
of GMRES iterations, indicating improved convergence rates. For example, for
\(CG(k)-DG(k)\) with \(k=2\), the average number of iterations decreases from
13.73 with \(n_{\text{smooth}}=1\) to 8.78 with \(n_{\text{smooth}}=2\). The
throughput shows a significant drop with an increase in the number of smoothing
steps. The decrease in the number of iterations does not balance the added
computational cost per iteration. Higher polynomial degrees result in lower
throughput due to the increased complexity of the smoother. The \(CG(k)-DG(k)\)
discretization generally requires more GMRES iterations compared to the
\(CG(k)-CGP(k)\) discretization for the same polynomial degree and smoothing
steps. The throughput for \(CG(k)-CGP(k)\) discretizations is consistently
higher than for \(CG(k)-DG(k)\) discretizations, suggesting better efficiency of
the continuous time discretization. At the same time, considering the point
evaluations plotted in~Figure~\ref{fig:rough-goal-quants}, we see that the
\(CG(k)-DG(k)\) discretization shows better convergence behavior compared to the
\(CG(k)-CGP(k)\) approach, due to the inhomogeneous material. As evidenced by
the same figure, both discretizations converge to the same solution. As the
resolution of the heterogeneities improves, the benefit of discontinuous
discretization becomes less pronounced. The advantage of higher-order
discretizations persists. The $k=2,\,r=6$ case achieves the highest throughputs,
with wall clock times that are slightly higher than those with $k=3,\,r=5$
(cf.~Figure~\ref{fig:rough-strong-scale}). However,
Figure~\ref{fig:rough-goal-quants} shows that the $k=3,\,r=5$ discretization
performs at least as well as the the $k=2,\,r=6$ one.

The high cost of the smoother represents a limitation in terms of the throughput
in degrees of freedom per second. In the convergence tests, higher order
discretizations demonstrated superior accuracy per work, while they now seem to
be disadvantagous in terms of the time to solution. We stress that this is
solely due to the increasingly expensive smoother. Therefore, future work should
focus on developing more efficient smoothers. A straightforward improvement
could be replacing the inner direct solver. Despite the high costs of the
smoother, it is quite robust as indicated by the consistent iteration counts
observed in all the performed tests. We note that the appeal of the proposed
space-time smoother is its wide applicability in fluid
mechanics~\cite{ahmedAssessmentSolversSaddle2018,anselmannGeometricMultigridMethod2023},
solid mechanics~\cite{wobkerNumericalStudiesVankaType2009}, fluid-structure
interaction~\cite{failerParallelNewtonMultigrid2021} and dynamic
poroelasticity~\cite{anselmannEnergyefficientGMRESMultigrid2024}. This
robustness together with the excellent scalability will enable the efficient
solution of large-scale coupled problems in future work.

\begin{figure}[htbp]
  \centering
  \includegraphics[width=\linewidth]{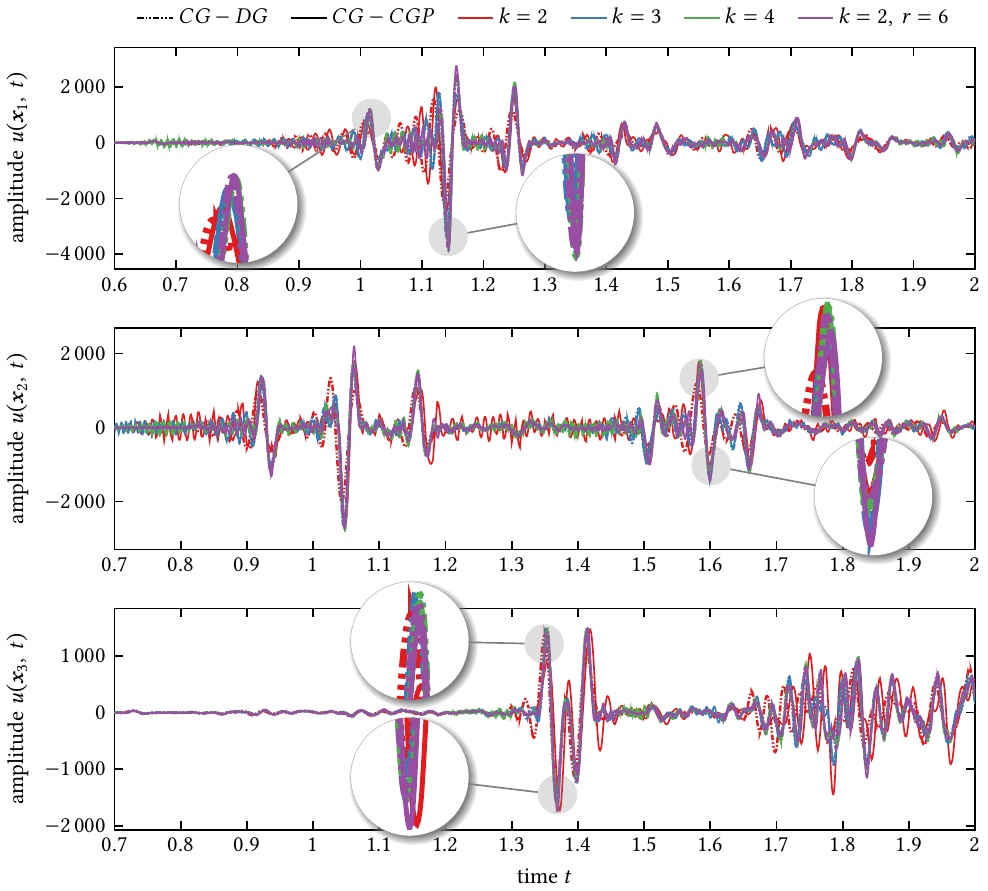}
  \caption{\label{fig:rougher-goal-quants}Plots of the point evaluations of the
    displacement at $\symbfit x_1,\:\symbfit x_2,\:\symbfit
    x_3$ plotted over subintervals of the time interval
    $[0,\,2]$. They are adapted such that the first signal arriving at the point
    is plotted.}
\end{figure}

\begin{figure}[htbp]
\includegraphics[width=\linewidth]{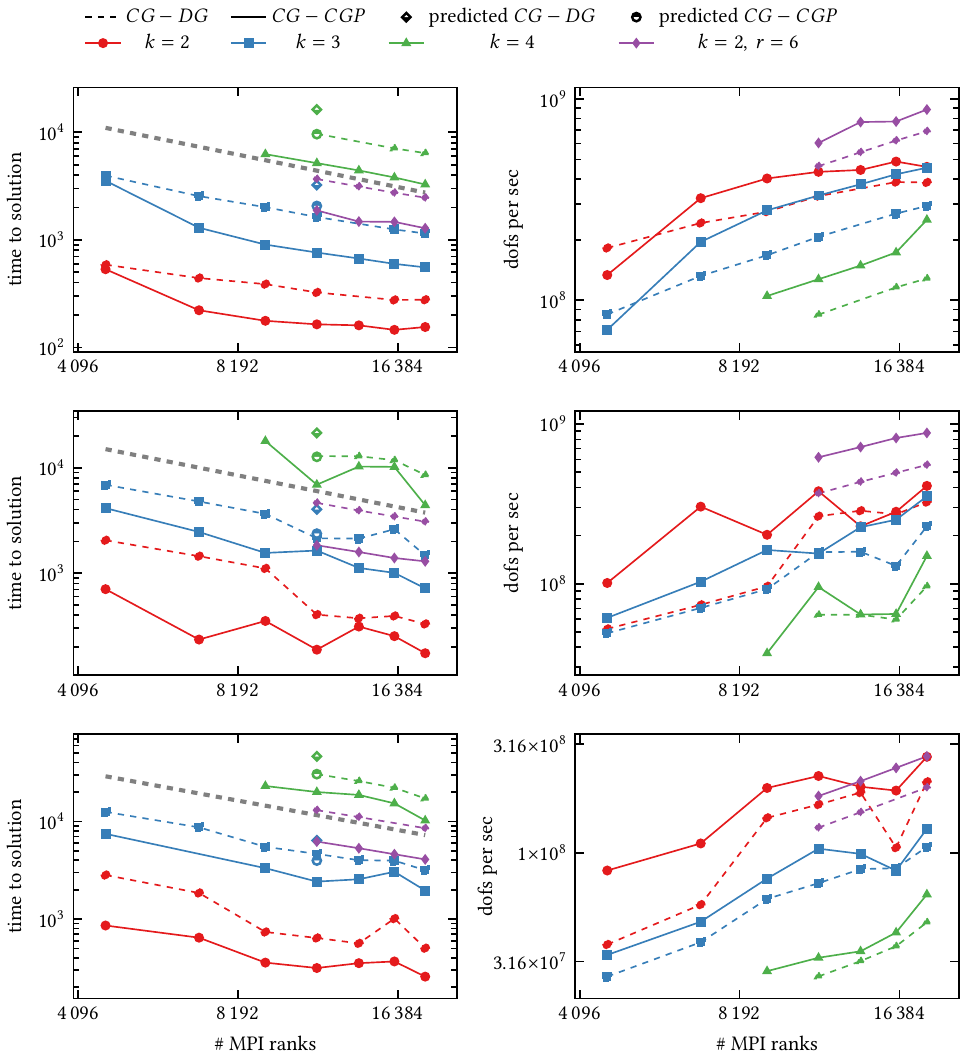}
\caption{\label{fig:rougher-strong-scale}Strong scaling test results for the
  STMG algorithm with varying numbers of smoothing steps. Each row corresponds
  to a different number of smoothing steps: 1, 2, and 4 from top to bottom. The
  left column shows the time to solution as a function of the number of MPI
  processes. The dashed gray lines indicate the optimal speedup. The
  half-squares and half-circles indicate the predicted runtime of the
  \(CG(k)\)-\(DG(k)\) and \(CG(k)\)-\(CGP(k)\) discretization, based on the
  wallclock time of one degree lower and~\eqref{eq:asm-complex}. The right
  column depicts the degrees of freedom (dofs) processed per second over the
  number of MPI processes.}
\end{figure}

\begin{figure}[htbp]
  \vspace*{-1ex}
  \centering%
  \includegraphics[width=.925\linewidth]{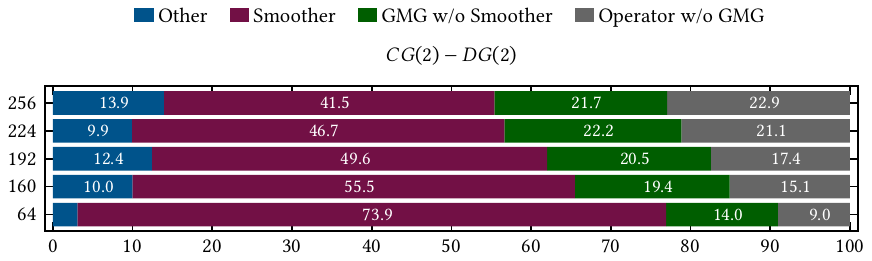}
  \includegraphics[width=.925\linewidth]{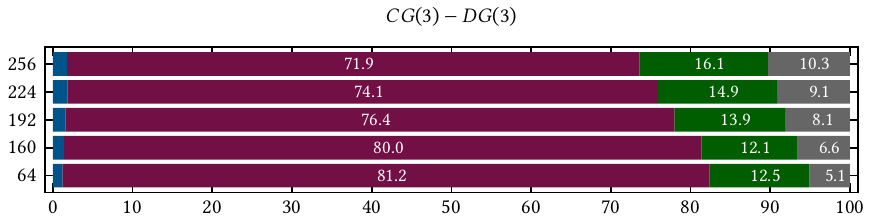}
  \includegraphics[width=.925\linewidth]{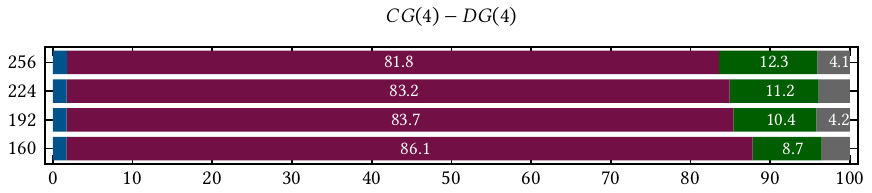}
  \includegraphics[width=.925\linewidth]{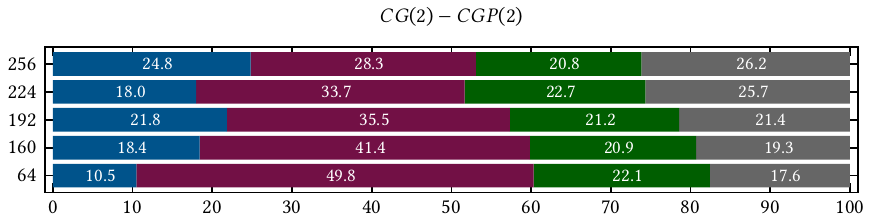}
  \includegraphics[width=.925\linewidth]{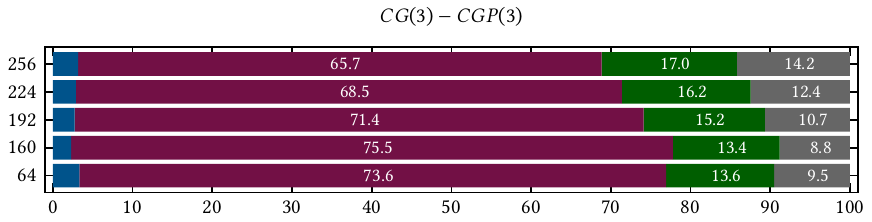}
  \includegraphics[width=.925\linewidth]{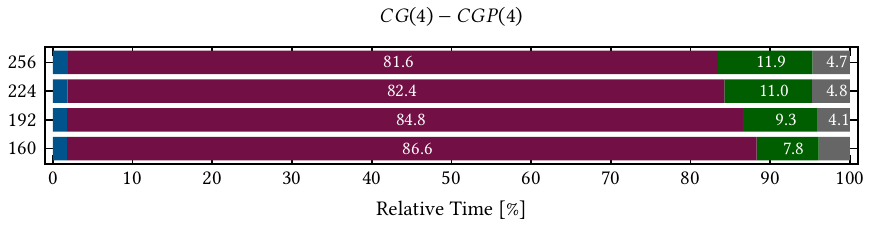}
\caption{\label{fig:rougher-strong-scale-rel}Time spent on different sections of
  the program in the strong scaling tests.}
\end{figure}

\FloatBarrier%
\paragraph*{Highly Heterogeneous Coefficients}
To increase the physical realism of our simulations, we incorporate highly
heterogeneous materials with rough coefficients and assess the method's ability
to effectively handle discontinuities and complex physical phenomena. To this
end, we perturb the coefficient function~\eqref{eq:rho} by sampling $c\in
[0.4,1.6]$ from a uniform random distribution for each coarse mesh cell and
using
$c\cdot\rho$ as the new coefficient. For the numerical examples we use the same
settings as for the heterogeneous coefficients. In
Figure~\ref{fig:rougher-strong-scale} we show the results of a scaling test
analogous to the one in Figure~\ref{fig:rough-strong-scale}. The number of
degrees of freedom, number of GMRES iterations and maximal throughput for this
setting are collected in Table~\ref{tab:iter-wave-rough-sss}. Analogous to the
previous test, we observe in~Figure~\ref{fig:rougher-goal-quants} that both
discretizations converge to the same solution. In
Figure~\ref{fig:rougher-strong-scale-rel} we present the relative timings of
different sections of the program during a strong scaling test with a single
smoothing step
$n_{\text{smooth}}=1$. Each bar represents a different run, and the segments of
the bar indicate the proportion of total execution time spent in various
components of the program. The ``GMG w/o Smoother'' section represents the time
spent in the STMG method, excluding the smoothing steps, i.\,e.\ only operator
evaluations and transfers. The proportion of this component remains almost
constant as the number of nodes increases. The ``Smoother'' section shows the
time spent in the ASM smoother. This is the computationally most intensive part
of the program, and the polynomial order has a significant impact on the time.
The proportion decreases slightly with an increasing number of nodes, but it is
always the dominant factor in wall clock time. The ``Operator w/o GMG'' section
shows the time spent in the operator, excluding the ones in the STMG
preconditioner. For lower order discretizations, this part becomes more
significant with an increasing number of nodes due to the communication in
combination with the low arithmetic intensity. The ``Other'' section represents
the remaining time not accounted for by the other components, i.\,e.\, the
assembly of source terms, evaluation of goal quantities, and operations
performed in between time steps. For lower order discretizations, this part also
increases with an increasing number of nodes due to the communication.

\begin{table}[htb]
  \centering
  \caption{Tables with the number of GMRES iterations until convergence for
    different polynomial degrees $k$ and number of refinements
    $r$. We show the numbers for the \(CG(k)-DG(k)\) and \(CG(k)-CGP(k)\)
    discretization of the wave equation. The last entry seperated by a rule is a value for $r=6$.}\label{tab:iter-wave-rough-sss}
    \begin{subcaptionblock}{\textwidth}
      \centering
      \caption{Iterations for \(CG(k)-DG(k)\) discretizations  with different numbers of smoothing steps ($n_{\text{smooth}}$).}\label{tab:iter-wave-dg-hhc-vss}
      \footnotesize
      \begin{tabular}{rrrrrrrrrr}
        \toprule
        $k$ & dofs$_{\Omega}$ &  dofs$_{\tau}$ &  dofs$_{\text{global}}$ & \multicolumn{2}{c}{$n_{\text{smooth}}=1$}                        & \multicolumn{2}{c}{$n_{\text{smooth}}=2$}                   & \multicolumn{2}{c}{$n_{\text{smooth}}=4$}\\
          &                &                &                       & $\overline{n}_{\text{iter}}$ & $\max \text{dofs}/\text{sec}$     & $\overline{n}_{\text{iter}}$& $\max \text{dofs}/\text{sec}$ & $\overline{n}_{\text{iter}}$& $\max \text{dofs}/\text{sec}$\\
        \midrule
        2 & \num{111284641} & $3\times 320$ & \num{106833255360} & 12.75 & \num{386656733} & 9.55  & \num{324031712}& 6.59  & \num{212392158}\\
        3 & \num{263374721} & $4\times 320$ & \num{337119642880} & 12.83 & \num{294170718} & 9.50  & \num{228555690}& 7.38  & \num{106548560}\\
        4 & \num{513922401} & $5\times 320$ & \num{822275841600} & 17.96 & \num{128380303} & 13.51 & \num{96341633}& 10.94 & \num{47862389}\\
        \midrule
        2 & \num{887503681} & $3\times 640$ & \num{1704007067520} & 14.25 & \num{691000432} & 9.48  & \num{553608534}& 6.54  & \num{199977358}\\
        \bottomrule
        \end{tabular}
    \end{subcaptionblock}%

  \begin{subcaptionblock}{\textwidth}
    \centering
    \caption{Iterations for \(CG(k)-CGP(k)\) discretizations with different numbers of smoothing steps.}\label{tab:iter-wave-cgp-hhc}
    \footnotesize
    \begin{tabular}{rrrrrrrrrr}
      \toprule
        $k$ & dofs$_{\Omega}$ &  dofs$_{\tau}$ &  dofs$_{\text{global}}$ & \multicolumn{2}{c}{$n_{\text{smooth}}=1$}                        & \multicolumn{2}{c}{$n_{\text{smooth}}=2$}                   & \multicolumn{2}{c}{$n_{\text{smooth}}=4$}\\
          &                &                &                       & $\overline{n}_{\text{iter}}$ & $\max \text{dofs}/\text{sec}$     & $\overline{n}_{\text{iter}}$& $\max \text{dofs}/\text{sec}$ & $\overline{n}_{\text{iter}}$& $\max \text{dofs}/\text{sec}$\\
      \midrule
      2 & \num{111284641} & $2\times 320$ &\num{71222170240} & 10.91 & \num{488827524} & 7.75  &\num{409322817}& 5.63 &\num{276805947} \\
      3 & \num{263374721} & $3\times 320$ &\num{252839732160}& 10.83 & \num{455567085} & 7.66  &\num{352291671} & 5.64 &\num{128999863} \\
      4 & \num{513922401} & $4\times 320$ &\num{657820673280}& 15.06 & \num{251168624} & 10.70 &\num{149165685} & 8.34 &\num{64366015} \\
      \midrule
      2 & \num{887503681} & $2\times 640$ & \num{1136004711680} &  11.65 & \num{886119119} & 7.80  & \num{877901632}& 5.66  & \num{278227948}\\
      \bottomrule
    \end{tabular}
  \end{subcaptionblock}%
\end{table}

The results of the tests for a problem with higher heterogeneity are summarized
in Table~\ref{tab:iter-wave-rough-sss}. The tests include \(CG(k)-DG(k)\) and
\(CG(k)-CGP(k)\) discretizations. Similar to our previous findings, we observe
that the number of GMRES iterations required for convergence increases with
higher polynomial degrees \(k\). Increasing the number of smoothing steps
\(n_{\text{smooth}}\) reduces the number of GMRES iterations. For example, for
\(CG(k)-DG(k)\) with \(k=2\), the average number of iterations decreases from
12.75 with \(n_{\text{smooth}}=1\) to 9.55 with \(n_{\text{smooth}}=2\). Again,
the \(CG(k)-DG(k)\) discretization generally requires more GMRES iterations
compared to the \(CG(k)-CGP(k)\) discretization for the same polynomial degree
and smoothing steps. The throughput for \(CG(k)-CGP(k)\) discretizations is
consistently higher than for \(CG(k)-DG(k)\) discretizations.

While more smoothing steps improve convergence by reducing the iteration count,
less smoothing steps are again beneficial in terms of the time to solution.
Overall, the results with larger inhomogeneity are consistent with our previous
findings. The larger inhomogeneity does not substantially change the convengence
and scalability of the method.

\FloatBarrier%
\section{Conclusion}
We verified the accuracy and robustness of the \(CG(p)-DG(k)\) and
\(CG(p)-CGP(k)\) methods for various polynomial orders on both Cartesian and
perturbed meshes. The proposed STMG method achieved optimal convergence and grid
independence. While we were not able to achieve independence of the polynomial
order, the increase with the number of iterations was moderate. We note that
robustness w.\,r.\,t.\ the polynomial order can be obtained by vertex-patch
based
smoothers~\cite{pavarinoAdditiveSchwarzMethods1993,miraiMultilevelAlgebraicError2020,schberlAdditiveSchwarzPreconditioning2008}
or $p$-multigrid~\cite{fehnHybridMultigridMethods2020}, which is a
straightforward extension to the space-time cell wise ASM smoother. However, the
critical bottleneck for higher polynomial orders is the computationally
expensive smoother. Thus, the inner direct solver should be replaced by a more
computationally less expensive alternative before turning to patch based
smoothers. In light of the preceding remarks, it is important to emphasize that
the method demonstrates remarkable performance. It achieves throughputs of over
a billion degrees of freedom per second on problems with more than a trillion
degrees of freedom and exhibits optimal scaling in our tests. Notably, it
outperforms existing matrix-based implementations by orders of magnitude.

We observed that, while an increase of the number of smoothing steps improves
convergence rates, it also has a negative impact on the throughput due to higher
computational costs. The reduction of GMRES iterations does not result in a
reduction of the time to solution, as the increased number of smoothing steps
counteract this effect. The \(CG(p)-CGP(k)\) discretization shows better
efficiency compared to \(CG(p)-DG(k)\). These advantages cannot be conclusively
verified in the context of problems with discontinuous material coefficients. To
investigate this further, we will examine the problems with discontinuous
material coefficients presented here with $DG(p)-DG(k)$ discretizations and
first-order formulations in space and time.

Higher-order discretizations achieve greater accuracy compared to lower-order
ones, owing to their efficiency in using fewer resources to attain the same
level of precision. In the convergence tests, we verified this in terms of the
accuracy per work. In the scaling tests, higher-order discretizations appeared
to be less advantageous in terms of time to solution. We established that the
high cost of the smoother is the reason for
this. %
Despite its high cost, the smoother exhibited excellent performance across
multiple applications. Overall, the proposed matrix-free space-time finite
element method within a space-time multigrid framework is efficient and
scalable, making it a promising candidate for large-scale problems in fluid
mechanics, fluid-structure interaction, and dynamic poroelasticity.
\subsection*{Acknowledgement}
Computational resources (HPC cluster HSUper) have been provided by the project
hpc.bw, funded by dtec.bw - Digitalization and Technology Research Center of the
Bundeswehr. dtec.bw is funded by the European Union - NextGenerationEU.

\printbibliography
\end{document}